\documentclass[a4paper,fleqn]{cas-sc}

\usepackage[numbers,sort&compress]{natbib}
\usepackage{amsmath,amssymb,amsfonts,amsthm}
\usepackage[fleqn]{subeqnarray}
\usepackage{graphicx}
\usepackage{tabularx}
\usepackage{subcaption}
\usepackage{algorithm}
\usepackage{algpseudocode}
\usepackage{enumitem}
\usepackage{float}
\usepackage{flafter}
\usepackage{placeins}
\usepackage{xcolor}
\usepackage{setspace}

\DeclareMathOperator*{\argmin}{arg\,min}

\begin{document}
\singlespacing
\let\WriteBookmarks\relax
\def\floatpagepagefraction{1}
\def\textpagefraction{.001}

\shorttitle{Residual-Based Physics-Aware MDNF for Inverse Problems}
\shortauthors{Lu et~al.}
\title[mode=title]{ResiPhy-MDNF: A Residual-Based Physics-Aware Multilevel Discrete Neural Field Framework for PDE-Constrained Inverse Problems}

\author[1]{Zheng Lu}
\ead{luzheng21@mails.jlu.edu.cn}

\author[1,3]{Jiwei Jia}
\cormark[1]
\ead{jiajiwei@jlu.edu.cn}

\author[2]{Young Ju Lee}
\cormark[2]
\ead{yjlee@txstate.edu}

\affiliation[1]{
  organization={School of Mathematics, Jilin University},
  addressline={No. 2699 Qianjin Street},
  city={Changchun},
  state={Jilin},
  postcode={130012},
  country={China}}

\affiliation[2]{
  organization={Department of Mathematics, Texas State University},
  addressline={MCS 470, 601 University Drive},
  city={San Marcos},
  state={TX},
  postcode={78666-4684},
  country={USA}}

\affiliation[3]{
  organization={Center for AI for Science and Engineering, Shenzhen Loop Area Institute},
  addressline={No. 6 Hongmian Road, Futian Free Trade Zone},
  city={Shenzhen},
  state={Guangdong},
  postcode={518038},
  country={China}}

\cortext[cor1]{Corresponding author.}
\cortext[cor2]{Corresponding author. This author was supported in part by NSF-DMS 2208499.}

\begin{abstract}
Inverse problems governed by partial differential equations are difficult when observations are sparse and the unknown coefficient field contains both large- and small-scale structures. We introduce a residual-based, physics-aware multilevel discrete neural field framework, ResiPhy-MDNF, for such problems. The method couples a coarse-to-fine discrete neural field (DNF) optimizer with a residual-based graph neural network (GNN) transfer operator, called ResiPhy-GNN. At each level, the DNF directly optimizes trainable grid- or mesh-based state and coefficient arrays using the prescribed numerical model. Between levels, ResiPhy-GNN maps the coarse representation to the fine representation through a learned prolongation based on graph connectivity, spatial features, and residual information. The method requires neither surrogate models nor offline pretraining. We evaluate the framework on coefficient inversion in Darcy flow for subsurface modeling and on electrical impedance tomography (EIT). In the controlled Darcy test case, the \(64^2\!\to128^2\) multilevel path uses \(1.25\times\) the cumulative grid-work proxy of the direct single-level \(128^2\) solve, while achieving \(6.76\times\) lower permeability error and \(10.4\times\) lower state error. On measured Kuopio Tomography Challenge 2023 EIT data, the framework improves the mean intersection-over-union after Otsu thresholding by roughly \(3.4\%\) over the official linearized complete electrode model reconstruction and by \(16.9\%\) over direct single-level discrete-field optimization. These results indicate that the same multilevel construction can be used across different coefficient structures, discretizations, and observation geometries.
\end{abstract}

\begin{highlights}
\item Discrete Neural Field (DNF) jointly optimizes discrete state and coefficient fields under fixed physics.
\item ResiPhy-GNN performs objective-fitted coarse-to-fine prolongation and correction.
\item Nested DNF levels provide a multilevel framework of coarse-to-fine inverse optimization.
\end{highlights}


\maketitle

\section{Introduction}
\label{sec:introduction}

Inverse problems governed by partial differential equations (PDEs) arise in
subsurface characterization, flow reconstruction, medical and industrial
imaging, material identification, and many other areas of computational
science.  In this work, we consider coefficient inverse problems concerned 
with recovering PDE coefficients from indirect and typically sparse observations 
of the corresponding solution fields.  Such problems are often ill
posed because the observation map is smoothing, the data dimension is small
relative to the number of coefficient unknowns, and distinct coefficient
fields may produce nearly indistinguishable observations
\cite{vogel2002computational,tarantola2005inverse,kaipio2005statistical}.
Their numerical solution combines physical modeling, regularization, and optimization algorithms.  These
difficulties are amplified on fine discretizations.  Although high-resolution 
discretizations are required to resolve interfaces,
inclusions, and channelized structures, they also introduce an increasing
number of weakly observable coefficient modes.
In a full-space formulation, the state and coefficient variables are optimized
jointly. Their strong coupling may allow different state--coefficient pairs to
produce similar observations, resulting in a poorly conditioned optimization
problem that is sensitive to initialization.

Finite-difference and finite-element inverse methods retain a prescribed
discrete forward model, but optimizing all fine-grid inverse variables from the
first iteration can be expensive.  Physics-informed neural networks (PINNs)
instead represent states and coefficients by coordinate networks and fit their
parameters through observation and residual losses
\cite{raissi2019physics,karniadakis2021physics}.  This representation may
suffer from spectral bias, imbalance among loss terms, and a mismatch between
continuous automatic differentiation and the discrete operator used to
interpret the measurements \cite{wang2022pinnfail}.  Learned surrogates and
neural operators reduce the cost of repeated forward evaluation, but they
usually require an offline dataset spanning the relevant coefficients,
sources, and observations \cite{lu2021deeponet,li2021fno}.  In coefficient
inversion, an error in the surrogate may be absorbed into the reconstructed
coefficient.  We therefore seek a method that keeps the numerical operator in
the inverse loss, uses no offline training data, and introduces fine-scale
unknowns progressively.

The last requirement is related to classical multilevel methods.  Full
multigrid initializes successive fine-grid solves with coarse-grid solutions
\cite{brandt2011fmg}, while algebraic multigrid constructs coarse spaces and
interlevel operators from a discrete algebraic system \cite{ruge1987amg}.
Subspace-correction and cascadic methods likewise treat error components on
the levels where they can be represented efficiently
\cite{xu1992subspace,bornemann1996cascadic}.  Multigrid ideas have also been
applied to discretized optimization and distributed parameter estimation
\cite{nash2000multigridoptimization,ascher2003multigridparameter,
borzi2009multigridoptimization}, and recent multilevel analyses address inverse
problems more generally \cite{weissmann2022multilevelinverse}.  Separately,
graph neural networks (GNNs) have been used to learn AMG prolongation and
multilevel domain-decomposition components
\cite{luz2020learningamg,taghibakhshi2023mggnn}.  These methods motivate a
coarse-to-fine inverse path, but they do not directly give an interlevel
correction fitted to the observations and discrete residual of a single
inverse instance.

In this paper, we introduce ResiPhy-MDNF.  Let
\(V_0\subset V_1\subset\cdots\subset V_J\) be nested reconstruction spaces associated with a hierarchy of increasingly fine discretizations (e.g., \(16^2\to32^2\to64^2\to128^2\) for Darcy or a coarse grid to an FEM mesh for EIT); a coarse-space function is naturally embedded in every finer space through prolongation \cite{brandt2011fmg}. On each level, a Discrete Neural Field (DNF)
solve represents state and unknown coefficient fields directly as trainable grid- or mesh-based
arrays and minimizes a prescribed finite-dimensional loss. The within-level 
DNF formulation is inspired by Optimizing a Discrete Loss
(ODIL) \cite{karnakov2024odil} and physics-informed cell-based representations 
\cite{kang2023pixel}, which directly optimizes discrete field
variables under a prescribed numerical objective. Successive DNF solves are 
coupled through ResiPhy-GNN, a residual-based,
physics-aware graph neural network built on message-passing operations
\cite{gilmer2017messagepassing}. At each transition from \(V_\ell\) to
\(V_{\ell+1}\), ResiPhy-GNN initializes the local prolongation weights
using the prescribed interpolation stencil, adjusts them through learned
edge-logit biases, and predicts a nodewise residual correction. The learned
weight adjustments and residual corrections are optimized jointly under the
target-level objective combining physics and data terms.

Consequently, our method can be applied to a range of PDE-constrained inverse
problems.

The main contributions are as follows.
\begin{enumerate}[leftmargin=2em]
    \item We formulate DNF optimization directly over the discrete inverse
    variables.  In the full-space Darcy realization, the pressure fields
    \(\{U_m\}_{m=1}^{M}\) and their shared permeability field \(K\) are
    optimized in one loss.  In the reduced Kuopio Tomography Challenge 2023
    (KTC2023) electrical impedance tomography (EIT) realization
    \cite{rasanen2024ktc2023,ktc2023zenodo}, the electrical states are
    eliminated by the linearized complete electrode model (CEM)
    \cite{somersalo1992existence}, and the conductivity contrast is optimized
    on the finite-element mesh.

    \item We construct ResiPhy-GNN as an objective-fitted, 
    learned prolongation operator that jointly adapts stencil weights and a nodewise residual update.

    \item We use nested DNF levels to provide a multilevel framework of
    coarse-to-fine inverse optimization.  Each level
    solves the inverse problem under its prescribed discrete loss, and the
    corrected interlevel representation initializes the next refined DNF solve.
\end{enumerate}

The experiments show three results.  First, on the
manufactured Darcy problem, direct \(128^2\) optimization and the
\(64^2\!\to128^2\) path both use 6000 final-grid Adam updates
\cite{kingma2015adam}.  The \(64^2\!\to128^2\) path uses a cumulative
grid-work proxy of \(1.25\), normalized by \(1.00\) for the direct
\(128^2\) solve, while achieving \(6.76\times\) lower coefficient error and
\(10.4\times\) lower state error.

Second, a matched \(64^2\!\to128^2\) LUNDIsim experiment compares simple
interpolation, learned weights only, and the complete ResiPhy-GNN transfer.
After the same number of refined DNF updates, ResiPhy-GNN reduces the final
state and permeability errors by approximately \(9\%\) and \(8\%\),
respectively.  The learned-weights-only ablation reduces the final state error
but does not improve the final permeability error.  A separate
\(16^2\!\to32^2\!\to64^2\!\to128^2\) experiment verifies the complete
multilevel path on a LUNDIsim dataset.

Third, we apply the method to KTC2023 EIT dataset.  The resulting
pipeline reaches \(0.623\) mean intersection over union after Otsu thresholding
\cite{otsu1979threshold}, compared with \(0.603\) for the official linearized
CEM reconstruction and \(0.533\) for direct single-level DNF.  This comparison 
evaluates the complete KTC2023 reconstruction pipeline; the matched LUNDIsim 
experiment provides the transfer ablation.  Taken together, the two
realizations show that the framework can be applied to different
coefficient structures, discrete models, and observation geometries.

The rest of the paper is organized as follows.
Section~\ref{sec:preliminaries} introduces the general formulations of
PDE-constrained inversion together with the Darcy and KTC2023 CEM-EIT models
and their discretizations. Section~\ref{sec:dnf} defines DNF, the nested
multilevel construction, and ResiPhy-GNN transfer.
Section~\ref{sec:problem_realizations} specifies the problem-dependent
trainable representations, reconstruction losses, and transfer configurations
for the Darcy and KTC2023 EIT realizations.
Section~\ref{sec:experiments} reports the numerical experiments, and
Section~\ref{sec:conclusion} summarizes the findings and limitations.
\section{General Formulation of PDE-Constrained Inversion}
\label{sec:preliminaries}

This section states the continuous and discrete formulations used below.  We
first introduce a PDE-constrained coefficient inverse problem, distinguish its
full-space, reduced-space, and linearized discrete forms, and then specify the
Darcy and electrical impedance tomography (EIT) models used in the two
realizations.

\subsection{PDE-constrained inverse problems}

Let \(\Omega\subset\mathbb{R}^{d}\) be a finite physical domain with boundary
\(\partial\Omega\), \(\kappa\) an unknown
coefficient field, and \(M\) the number of experiments.  For each experiment
\(m\in\{1,\ldots,M\}\), let \(s_m\) collect the known driving inputs (source terms, boundary values, or applied electrode currents) and
let \(u_m\) denote the corresponding state field.  A general forward and
observation model can be written as
\begin{subeqnarray}
    \mathcal{F}(u_m,\kappa;s_m)
    & = & 0 \qquad \text{in }\Omega,
        \slabel{eq:prelim_continuous_pde}\\
    \mathcal{B}(u_m,\kappa;s_m)
    & = & 0 \qquad \text{on }\partial\Omega,
        \slabel{eq:prelim_continuous_bc}\\
    d_m
    & = & \mathcal{P}_m u_m+\eta_m,
        \qquad m=1,\ldots,M.
        \slabel{eq:prelim_observation}
\end{subeqnarray}
Here \(d_m\) is the measured data for experiment
\(m\), with \(\mathcal P_m\) the observation operator and \(\eta_m\) the
measurement error.  \(\mathcal F\) and \(\mathcal B\) are fixed operators representing the governing differential equation and its boundary or auxiliary conditions.  The coefficient \(\kappa\) is
an unknown shared by all experiments, whereas each experiment has its own
state \(u_m\).  Through \(\mathcal P_m\), measurements may be taken on an
interior subset of \(\Omega\) or a boundary subset of \(\partial\Omega\).

A regularized full-space formulation retains the states and coefficient as
independent unknown variables and reads as follows
\cite{lions1971optimal,hinze2009optimization,borzi2012computational}:
\begin{equation}
    \begin{aligned}
    \min_{ \left \{ \kappa,\{u_m\}_{m=1}^{M} \right \} }\quad
    &
    \frac{1}{2}\sum_{m=1}^{M}
    \|\mathcal{P}_m u_m-d_m\|_{W_m}^{2}
    +\lambda_{\rm reg}\,\mathcal{R}(\kappa) \\
    \text{subject to}\quad
    &
    \mathcal{F}(u_m,\kappa;s_m)=0
    \quad \text{in }\Omega,\\
    &
    \mathcal{B}(u_m,\kappa;s_m)=0
    \quad \text{on }\partial\Omega,\qquad m=1,\ldots,M .
    \end{aligned}
    \label{eq:prelim_continuous_inverse}
\end{equation}
Here \(W_m\) is the positive-semidefinite data-weighting matrix, interpreted
statistically as a measurement precision matrix
\cite{kaipio2005statistical}.  We define
\(\|v\|_{W_m}^2:=v^T W_m v\), which is a seminorm if \(W_m\) is singular;
\(\lambda_{\rm reg}\geq0\) is the regularization weight, and
\(\mathcal{R}\) encodes prior information such as smoothness, total-variation
structure, or proximity to a background field.  Regularization is needed
because elliptic coefficient inverse problems are typically ill posed:
different fine-scale coefficients may have very similar effects on sparse
interior or boundary measurements.

\subsection{Discrete full- and reduced-space formulations}

Classical numerical inversion discretizes the forward problem and then solves
the resulting finite-dimensional optimization problem
\cite{lions1971optimal,hinze2009optimization,borzi2012computational}.  Let
\(h\) index the discretization, let \(n_\kappa\) and \(n_u\) denote the numbers
of coefficient and state degrees of freedom, and let
\(\kappa_h\in\mathbb{R}^{n_\kappa}\) and
\(U_{m,h}\in\mathbb{R}^{n_u}\) denote the discrete coefficient and state,
respectively.  The discrete forward problem is written as
\begin{equation}
    F_h(U_{m,h},\kappa_h;s_{m,h})=0,\qquad m=1,\ldots,M,
    \label{eq:prelim_discrete_constraint}
\end{equation}
where \(F_h\) is the fixed algebraic residual operator, \(s_{m,h}\) is the
discrete excitation data, and \(P_{m,h}\) maps the state degrees of freedom to
the data space.  The same data-space weighting matrix \(W_m\) is retained
after discretization.  The discrete full-space formulation is
\begin{equation}
    \begin{aligned}
    \min_{\kappa_h,\{U_{m,h}\}_{m=1}^{M}}\quad
    &
    \frac{1}{2}\sum_{m=1}^{M}
    \|P_{m,h}U_{m,h}-d_m\|_{W_m}^{2}
    +\lambda_{\rm reg} R_h(\kappa_h)\\
    \text{s.t.}\quad
    &
    F_h(U_{m,h},\kappa_h;s_{m,h})=0,\qquad m=1,\ldots,M.
    \end{aligned}
    \label{eq:prelim_full_space}
\end{equation}

Here \(R_h\) denotes the discrete regularizer. In the full-space formulation
\eqref{eq:prelim_full_space}, both the states
\(\{U_{m,h}\}_{m=1}^{M}\) and the shared coefficient \(\kappa_h\) are treated
as optimization variables.

For a measurement map linearized around a background coefficient, let
\(\delta\kappa_h\) denote a coefficient perturbation, \(\delta d\) the
corresponding stacked data perturbation, and \(J_h\) the discrete sensitivity
or Jacobian of the measurement map.  Then
\begin{equation}
    \delta d\approx J_h\delta\kappa_h.
    \label{eq:prelim_linearized_map}
\end{equation}
A standard Tikhonov update
\cite{vogel2002computational,kaipio2005statistical} is defined by
\begin{equation}
    \min_{\delta\kappa_h}
    \frac{1}{2}\|J_h\delta\kappa_h-\delta d\|_{W}^{2}
    +\frac{\lambda}{2}\|L_h\delta\kappa_h\|_2^{2},
    \label{eq:prelim_tikhonov}
\end{equation}
where \(W\) is the stacked counterpart of \(W_m\), \(L_h\) is the
regularization operator,
\(\lambda\geq0\) is its weight, and \(\|\cdot\|_2\) is the Euclidean norm.

The full-space, reduced-space, and linearized formulations differ in the
variables retained by the optimization, but each is defined by a prescribed
numerical operator, observation map, and regularization.  The following two
models instantiate these formulations: Darcy inversion is posed in full space,
whereas complete-electrode-model EIT is posed as a reduced linearized inverse
problem.

\subsection{Darcy model and discretization}
\label{subsec:darcy_model_discretization}

For permeability \(K(x)>0\), source term \(f_m\), and prescribed boundary
pressure \(g_m\), the pressure \(u_m\) generated by excitation \(m\) satisfies
\begin{equation}
    -\nabla\cdot(K\nabla u_m)=f_m
    \quad\text{in }\Omega,
    \qquad
    u_m=g_m
    \quad\text{on }\partial\Omega .
    \label{eq:prelim_darcy_model}
\end{equation}
On a structured grid, let \(C_i\) be cell \(i\), \(|C_i|\) its volume,
\(\mathcal{N}(i)\) the indices of its face-neighboring cells, and
\(U_{m,h,i}\), \(K_{h,i}\), and \(f_{m,h,i}\) the cell values of pressure,
permeability, and source, respectively.  Integration over \(C_i\) gives the
conservative residual
\begin{equation}
    \operatorname{Res}^{\rm D}_{m,h,i}(U_{m,h},K_h)
    =
    \sum_{j\in\mathcal{N}(i)}
    T_{ij}(K_h)(U_{m,h,i}-U_{m,h,j})
    -|C_i|f_{m,h,i},
    \label{eq:prelim_darcy_discrete}
\end{equation}
Here the superscript \({\rm D}\) denotes Darcy flow, while \(m\), \(h\), and
\(i\) identify the experiment, discretization, and control volume \(C_i\),
respectively.  The prescribed boundary contributions are included in the same balance.
For an orthogonal grid, the interior transmissibility is
\[
    T_{ij}(K_h)
    =
    \frac{|e_{ij}|}{d_{ij}}
    \frac{2K_{h,i}K_{h,j}}{K_{h,i}+K_{h,j}},
\]
where \(e_{ij}\) is the common face of cells \(i\) and \(j\), \(|e_{ij}|\) is
its measure, and \(d_{ij}\) is the distance between the cell centers.  The
harmonic face value preserves flux continuity across permeability contrasts.
For each source or boundary excitation \(m=1,\ldots,M\), sparse pressure
observations have the form
\(d_m=P_{m,h}U_{m,h}+\eta_m\).

\subsection{KTC2023 CEM-EIT model and discretization}
\label{subsec:eit_model_discretization}

Let \(\sigma(x)>0\) denote conductivity and let \(u_m\) be the interior
electric potential generated by current pattern \(m\).  Let \(E\) be the
number of electrodes, \(e\in\{1,\ldots,E\}\) the electrode index,
\(E_e\subset\partial\Omega\) the electrode surfaces, and
\(\Gamma_{\rm el}:=\bigcup_{e=1}^{E}E_e\) the electrode-covered portion of
the boundary.  Let \(\zeta_e>0\) be the contact impedances,
\(I_{m,e}\) the applied currents, and \(V_{m,e}\) the electrode voltages.  The applied currents satisfy
\(\sum_{e=1}^{E}I_{m,e}=0\).  If \(n\) is the outward unit normal, then
\(\partial_n u_m:=\nabla u_m\cdot n\), and \(ds\) denotes the boundary
integration measure.  The complete electrode model (CEM)
\cite{somersalo1992existence} is
\begin{subequations}
\label{eq:prelim_cem}
\begin{alignat}{2}
    \nabla\cdot(\sigma\nabla u_m)
    &=0
    &\qquad& \text{in }\,\Omega,\\
    u_m+\zeta_e\sigma\partial_n u_m
    &=V_{m,e}
    && \text{on }\,E_e,\qquad e=1,\ldots,E,\\
    \int_{E_e}\sigma\partial_n u_m\,ds
    &=I_{m,e}
    && e=1,\ldots,E,\\
    \sigma\partial_n u_m
    &=0
    && \text{on }\,\partial\Omega\setminus\Gamma_{\rm el},\\
    \sum_{e=1}^{E}V_{m,e}
    &=0.&&
\end{alignat}
\end{subequations}
A finite-element discretization gives the grounded block system
\begin{equation}
    A_h(\sigma_h)
    \begin{bmatrix}U_{m,h}\\V_{m,h}\end{bmatrix}
    =b_m,
    \qquad
    y_{m,h}=Q_mV_{m,h},
    \label{eq:prelim_cem_discrete}
\end{equation}
where \(\sigma_h\) is the discrete conductivity,
\(A_h(\sigma_h)\) is the conductivity-dependent system matrix,
\(U_{m,h}\) and \(V_{m,h}\) collect the nodal and electrode potentials,
respectively, \(b_m\) encodes the applied current pattern, \(Q_m\) selects
the measured voltage channels, and \(y_{m,h}\) is the predicted measurement
vector.  For difference imaging about a homogeneous background conductivity
\(\sigma_0\), the stacked voltage change is linearized as
\begin{equation}
    \Delta V\approx J_h\delta\sigma_h,
    \qquad
    \delta\sigma_h=\sigma_h-\sigma_{0,h},
    \label{eq:prelim_cem_linearization}
\end{equation}
where \(\Delta V\) is the stacked measured electrode-voltage change,
\(\sigma_{0,h}\) is the discrete background conductivity,
\(\delta\sigma_h\) is the discrete conductivity contrast, and \(J_h\) is the
finite-element voltage-sensitivity matrix.  This is a reduced linear inverse model of the form
\eqref{eq:prelim_linearized_map}.

For KTC2023, the official active voltage mask is applied and the current
patterns are stacked.  We suppress the mesh subscript and write
\(J\delta\sigma\approx\Delta V\), where \(J\) is the official stacked
finite-element sensitivity matrix and \(\delta\sigma=\sigma-\sigma_0\) is the
conductivity contrast
\cite{kaipio2005statistical,rasanen2024ktc2023,ktc2023zenodo}.
These prescribed discretizations provide the application models used to
instantiate the framework developed in Section~\ref{sec:dnf}; their
problem-specific trainable representations, losses, and transfer features are
specified in Section~\ref{sec:problem_realizations}.
\section{The ResiPhy-MDNF Framework: DNF Optimization and ResiPhy-GNN Transfer}
\label{sec:dnf}
\emph{ResiPhy-MDNF} is a multilevel inverse framework that retains the prescribed numerical model and optimizes each inverse instance without offline training data. It has three components.
First, at each resolution level, the active inverse variables are represented directly as trainable 
grid- or mesh-based fields and optimized under the corresponding discrete loss function; we refer to this 
within-level procedure as a \emph{Discrete Neural Field} (DNF) solve. Second, these DNF solves are 
arranged over a nested coarse-to-fine hierarchy, forming the \emph{multilevel discrete neural field} (MDNF) 
backbone. Third, consecutive levels are connected by a single interlevel transfer operator that performs
interpolation-initialized learned prolongation with a residual update. This transfer operator is fitted using the
current observations and the target-level loss function, without fine-grid
reference fields. We call this transfer module the residual-based,
physics-aware graph neural network, abbreviated as ResiPhy-GNN
\cite{gilmer2017messagepassing}; it provides the \emph{ResiPhy} component in
the framework name.

Coarse-to-fine transfer is intrinsically local: each fine-level degree of freedom 
draws information from only a few coarse neighbors defined by the interpolation 
stencil. A graph neural network makes this locality explicit by aggregating messages 
from those neighbors and sharing the same update rule across target 
nodes \cite{gilmer2017messagepassing}. This design preserves 
the sparse connectivity of the transfer operator—unlike a dense MLP or a globally 
learned interpolation map—and naturally handles irregular meshes without requiring 
the translation-invariant Cartesian structure that typical CNNs assume \cite{pfaff2021learning}.

Each coarse-to-fine transition consists of a current-level DNF solve followed
by ResiPhy-GNN transfer. After the active inverse variables have been optimized
on the current discretization, ResiPhy-GNN constructs an incoming
representation on the target discretization, which initializes the next-level
DNF solver. The reconstruction and transfer objectives use the same numerical
operator, boundary treatment, and measurement model. ResiPhy-GNN supplies the
interlevel initialization, while the next-level DNF solve is governed by its
reconstruction objective. The corresponding
problem-specific representations, objectives, and features are presented in
Section~\ref{sec:problem_realizations}.

\paragraph{Notation.}
The framework uses five recurring symbols.  At level \(\ell\),
\(\Theta_\ell\) denotes the trainable variables, \(z_\ell\) the assembled
physical fields, \(\tau_\ell\) the transfer coordinates associated with
level \(\ell\), and \(\mathcal L_\ell\) the reconstruction loss.  At the
interface from level \(\ell\) to level \(\ell+1\),
\(\mathcal T_{\ell+1}\) denotes the transfer loss.  Unless marked by \({\rm opt}\), \(\Theta_\ell\) denotes the
trainable variables initialized at the beginning of level \(\ell\);
\(\Theta_\ell^{\rm opt}\) and \(z_\ell^{\rm opt}\) denote the optimized
variables and assembled fields.  The reconstruction spaces satisfy
\(V_0\subset\cdots\subset V_J\), from coarse to fine.  Three fixed,
application-dependent operations connect these objects:
\begin{equation}
    z_\ell=\operatorname{Asm}_\ell(\Theta_\ell),\qquad
    \tau_\ell=\operatorname{Out}_\ell(\Theta_\ell^{\rm opt}),\qquad
    \Theta_{\ell+1}=\operatorname{In}_{\ell+1}(\tau_{\ell+1}).
    \label{eq:framework_notation_roadmap}
\end{equation}
\(\operatorname{Asm}_\ell\) maps raw trainable variables to physically valid fields: it enforces parameter constraints (e.g., positivity via softplus), restores prescribed boundary values, and, when the trainable representation and the numerical discretization live on different grids, maps between the two. \(\operatorname{Out}_\ell\) extracts the variables used by the interlevel transfer, and \(\operatorname{In}_{\ell+1}\) converts the transferred representation into the next-level trainable variables. The concrete forms of these three maps are given in Section~\ref{sec:problem_realizations}.
Symbols local to an analysis or realization are defined where they occur.
\subsection{Framework overview}
\label{subsec:framework_overview}

On level \(\ell\), the assembled field \(z_\ell\) can contain both state variables
and coefficient variables.  In a reduced problem, it may contain only a
coefficient or conductivity perturbation after the state has been eliminated by
a forward model or a linearization.  The assembly operator \(\operatorname{Asm}_\ell\)
inserts prescribed
values, applies positivity or box constraints, maps reduced
parameters to the active reconstruction space, and returns the admissible
fields seen by the physics evaluator.

The transfer coordinates \(\tau_\ell\) may differ from the assembled field \(z_\ell\): for Darcy,
\(\tau_\ell\) contains raw \(\rho_\ell\) whereas \(z_\ell\) contains the positive \(K_\ell\); for EIT,
the transfer passes through an auxiliary grid before reaching the FEM mesh.
\(\operatorname{Out}_\ell\) extracts the variables acted on by transfer;
\(\operatorname{In}_{\ell+1}\) converts the corrected representation into trainable variables.
Their domains need not coincide.

Let \(d\) collect the observations.  On level \(\ell\), DNF minimizes the
following loss directly over \(\Theta_\ell\):
\begin{equation}
    \label{eq:mdnf_generic_loss}
    \mathcal{L}_\ell(\Theta_\ell)
    =\mathcal{J}^{\rm data}_\ell(z_\ell)
    +\lambda_{\rm pde}\,\mathcal{J}^{\rm phys}_\ell(z_\ell)
    +\lambda_{\rm reg}\,\mathcal{R}_\ell(z_\ell),
    \qquad z_\ell=\operatorname{Asm}_\ell(\Theta_\ell).
\end{equation}
Here the three terms measure data mismatch, physics inconsistency, and prior
regularity.  The nonnegative scalars
\(\lambda_{\rm pde}\) and \(\lambda_{\rm reg}\) weight the physics and regularization terms,
respectively.  The penalty formulation is a standard approach for PDE-constrained optimization \cite{nocedal2006numerical}.  Dependence on the fixed observations is omitted from the notation.  The
physics evaluator is part of this loss, not a learned module, and automatic
differentiation propagates through \(\operatorname{Asm}_\ell\).

Concretely, one coarse-to-fine transition maps the outgoing representation
$\tau_\ell$ (extracted from the completed DNF solve) to an incoming
representation $\tau_{\ell+1}$, which is an initialization associated with
$V_{\ell+1}$, not the optimized solution on that space.  The maps
$\operatorname{Out}_\ell$ and $\operatorname{In}_{\ell+1}$ connect the DNF
variables to this transfer; the ResiPhy-GNN weights $\omega_\ell$ are fitted
with a fixed transfer loss $\mathcal T_{\ell+1}$, distinct from the next
level's reconstruction loss.  The transition chain is:

\begin{equation}
    \Theta_\ell
    \xrightarrow{\;\mathrm{DNF}\;}
    (\Theta_\ell^{\rm opt},z_\ell^{\rm opt})
    \xrightarrow{\;\operatorname{Out}_\ell\;}
    \tau_\ell
    \xrightarrow{\;\mathrm{ResiPhy\text{-}GNN}\;}
    \tau_{\ell+1}
    \xrightarrow{\;\operatorname{In}_{\ell+1}\;}
    \Theta_{\ell+1} .
    \label{eq:complete_level_transition}
\end{equation}

Algorithm~\ref{alg:complete_mdnf_gnn} gives the complete order of operations.
Deterministic interpolation is recovered when the learned edge-logit biases
and residual outputs are both zero.

\begin{algorithm}[H]
\caption{ResiPhy-MDNF: single-level DNF optimization with ResiPhy-GNN transfer}
\label{alg:complete_mdnf_gnn}
\begin{algorithmic}[1]
    \Require Nested coefficient or contrast spaces
    \(V_0\subset\cdots\subset V_J\), data \(d\), fixed
    reconstruction losses \(\mathcal L_\ell\), fixed
    transfer losses \(\{\mathcal T_{\ell+1}\}_{\ell=0}^{J-1}\), assembly maps
    \(\operatorname{Asm}_\ell\), outgoing maps
    \(\operatorname{Out}_\ell\), incoming maps
    \(\operatorname{In}_\ell\), and the ResiPhy-GNN architecture
    \State Initialize trainable variables \(\Theta_0\) associated with \(V_0\)
    \For{\(\ell=0,\ldots,J\)}
        \State \textbf{DNF solve:} optimize \(\Theta_\ell\) with the fixed DNF
        loss \(\mathcal{L}_\ell\)
        \State Assemble
        \(z_\ell^{\rm opt}=\operatorname{Asm}_\ell(\Theta_\ell^{\rm opt})
        \in\mathcal{Z}_\ell\)
        \If{\(\ell=J\)}
            \State \Return \(z_J^{\rm opt}\in\mathcal{Z}_J\)
        \EndIf
        \State Extract transfer coordinates
        \(\tau_\ell=\operatorname{Out}_\ell(\Theta_\ell^{\rm opt})\)
        \State \textbf{ResiPhy-GNN transfer:} initialize the coarse-to-fine
        weights from interpolation and assemble nodewise stencil features
        \State Reinitialize and fit only \(\omega_\ell\) for the current
        inverse instance by minimizing \(\mathcal T_{\ell+1}\)
        evaluated at the ResiPhy-GNN output
        \State Set \(\tau_{\ell+1}\) to the fitted
        ResiPhy-GNN output generated from \(\tau_\ell\) and \(d\)
        \State Initialize
        \(\Theta_{\ell+1}=
        \operatorname{In}_{\ell+1}(\tau_{\ell+1})\)
    \EndFor
\end{algorithmic}
\end{algorithm}

\begin{figure}[!htbp]
\centering
\includegraphics[width=0.95\linewidth]{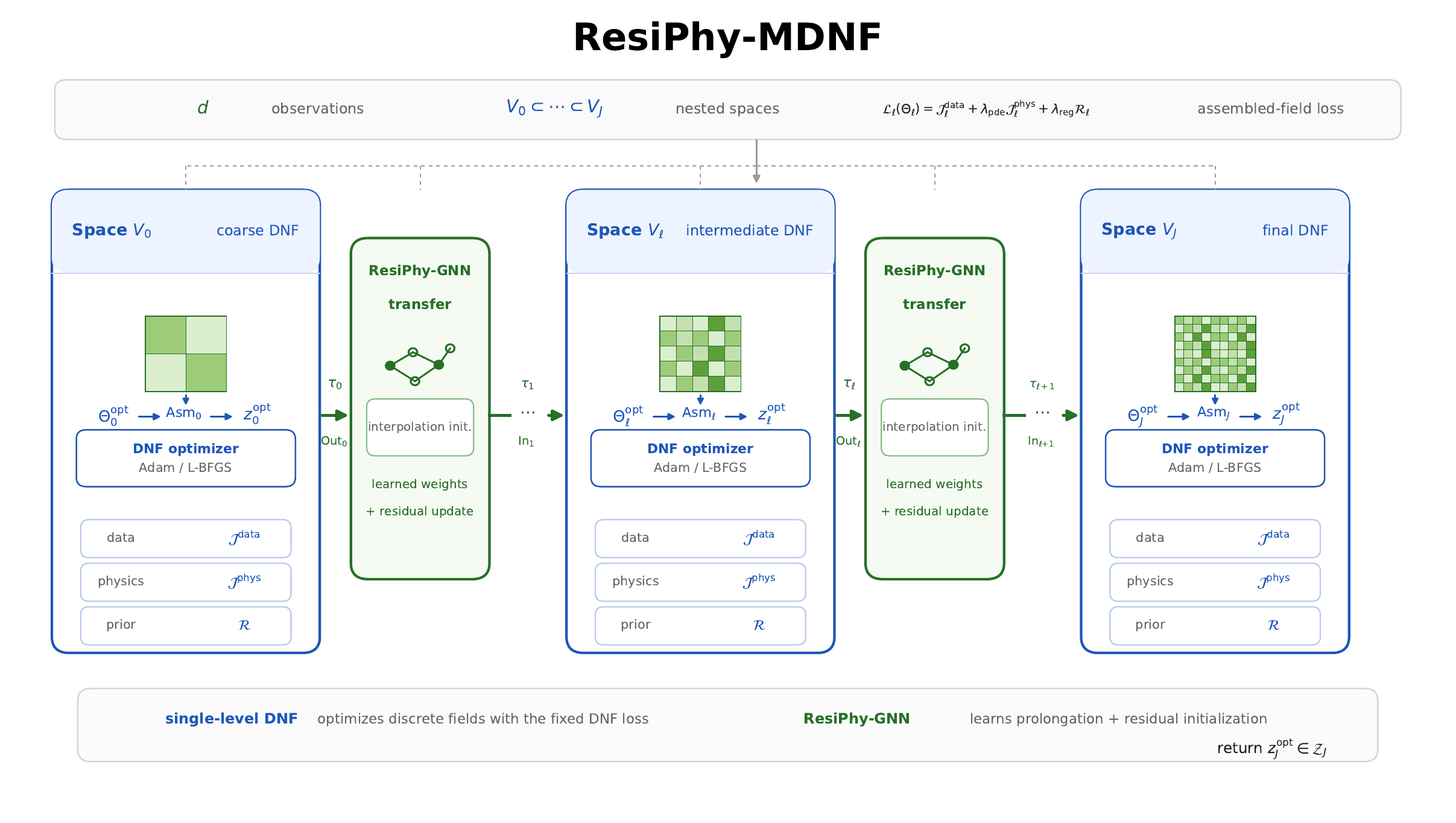}
\caption{Algorithmic view of ResiPhy-MDNF.  Each DNF level optimizes its
assembled fields, and ResiPhy-GNN maps the outgoing representation
\(\tau_\ell\) to an incoming representation \(\tau_{\ell+1}\) through
interpolation-centered learned prolongation and a residual update.}
\label{fig:dnf_architecture}
\end{figure}

\subsection{Single-level DNF optimization}
\label{subsec:one_level_dnf}

The full-space and reduced-space formulations in
Section~\ref{sec:preliminaries} distinguish which variables remain active on
one fixed discretization.  The DNF step in
Algorithm~\ref{alg:complete_mdnf_gnn} is the only step that
optimizes the inverse variables on a fixed level.  Starting from
\(\Theta_\ell\), it repeatedly assembles admissible fields, evaluates the fixed
DNF loss in \eqref{eq:mdnf_generic_loss}, differentiates through the numerical
operations, and updates \(\Theta_\ell\) with Adam
\cite{kingma2015adam} or the limited-memory
Broyden--Fletcher--Goldfarb--Shanno (L-BFGS) method
\cite{liu1989lbfgs}.  The output is the
optimized assembled field \(z_\ell^{\rm opt}\in\mathcal{Z}_\ell\), not a network
prediction.

\begin{algorithm}[H]
    \caption{DNF optimization on one level}
    \label{alg:dnf_single_level}
    \begin{algorithmic}[1]
        \Require Trainable variables \(\Theta_\ell\), assembly map
        \(\operatorname{Asm}_\ell\), loss terms in \eqref{eq:mdnf_generic_loss}
        \For{optimization step \(t=0,1,\ldots\)}
            \State Assemble admissible fields
            \(z_\ell=\operatorname{Asm}_\ell(\Theta_\ell)\in\mathcal{Z}_\ell\)
            \State Evaluate data, physics, and regularization terms
            \State Differentiate \(\mathcal{L}_\ell(\Theta_\ell)\) by automatic differentiation
            \State Update \(\Theta_\ell\) with Adam or L-BFGS
        \EndFor
        \State \Return optimized fields \(z_\ell^{\rm opt}\in\mathcal{Z}_\ell\)
    \end{algorithmic}
\end{algorithm}

The trainable objects are the discrete fields themselves.  In a
full-space realization, the assembled field contains states and coefficients;
in a reduced realization, it contains the parameters seen by a preassembled
measurement map.  The Darcy and EIT choices of
\(\mathcal{J}^{\rm phys}_\ell\) are specified in
Section~\ref{subsec:darcy_realization} and
Section~\ref{subsec:eit_realization}.

\subsection{Multilevel DNF backbone}
\label{subsec:mdnf_backbone}

MDNF uses a hierarchy of the active coefficient or contrast space
across multiple levels.
Elliptic coefficient
inverse problems are ill conditioned because fine-scale coefficient modes can
be weakly visible in sparse or boundary data.  MDNF therefore repeats the
single-level DNF solve over the nested coefficient or contrast reconstruction
spaces.
\[
    V_0\subset V_1\subset\cdots\subset V_J,
\]
rather than introducing all fine-scale degrees of freedom in the first solve.
Here \(V_\ell\) is the active coefficient or contrast reconstruction space at
level \(\ell\), and \(V_J\) is the final coefficient space.  Let
\(n_\ell=\dim V_\ell\), and let
\(\vartheta_\ell\in\mathbb{R}^{n_\ell}\) denote the trainable coordinates of
the coefficient or contrast on that level.  If state variables are active,
they are included together with \(\vartheta_\ell\) in \(\Theta_\ell\).  The
assembly operator introduced above converts these trainable variables into the
physical fields:
\[
    z_\ell=\operatorname{Asm}_\ell(\Theta_\ell)
    \in\mathcal{Z}_\ell.
\]
The application-specific coefficient or contrast assembly is given in
Section~\ref{sec:problem_realizations}.  For Darcy,
\(\vartheta_\ell=\rho_\ell\) and
\(K_\ell=K_{\min}+\operatorname{softplus}(\rho_\ell)\); for EIT,
\(\vartheta_\ell\) represents a signed conductivity contrast that is placed on
the prescribed finite-element space by a fixed grid-to-mesh map.  The loss
restricted to the current level is the DNF loss
\(\mathcal{L}_\ell(\Theta_\ell)\).

\subsection{Residual-based physics-aware GNN transfer}
\label{subsec:residual_gnn_framework}

The MDNF backbone in Section~\ref{subsec:mdnf_backbone} defines the nested
spaces and the rationale for coarse-to-fine continuation.
To complete the hierarchy, however, an interlevel
operator is still needed to convert a completed DNF solution into the incoming
representation for the next level.  This operator combines
interpolation-initialized, learned coarse-to-fine aggregation with a residual
GNN transfer map, where GNN
denotes a graph neural network.  This GNN-based transfer module is called the
residual-based, physics-aware graph neural network (ResiPhy-GNN).  It operates between two DNF levels.

\paragraph{ResiPhy-GNN architecture.}
Interlevel transfer is a local operation: each fine node depends on a small 
coarse stencil via prescribed interpolation weights $w^P_{i_fi_c}$. Fixed 
interpolation, however, ignores the observations and PDE residual available at 
the target level. ResiPhy-GNN therefore uses a shared pointwise network 
$\operatorname{Corr}_{\omega_\ell}$ to provide instance-aware corrections
through two mechanisms: learned edge-logit biases $b_{i_fi_c}$ that adjust 
the aggregation weights from their geometric baseline, and a node residual 
$\operatorname{res}_{i_f}$ that adds a bounded additive correction. The 
interpolation-centered softmax $\widehat w_{i_fi_c}\propto w^P_{i_fi_c}\exp(b_{i_fi_c})$ 
ensures that when $b_{i_fi_c}=0$ the transfer recovers deterministic 
interpolation exactly; the network therefore needs only to learn the 
physics-dependent correction relative to a known valid baseline. The shared 
pointwise application enforces translation invariance—the correction depends 
only on the local node context—and makes the number of trainable parameters 
independent of discretization size.

At interface \(\ell\), let \(i_c\) and \(i_f\) denote coarse-node and
fine-node indices, respectively.  Let \(w^P_{i_fi_c}\) denote the prescribed
interpolation weight from coarse node \(i_c\) to fine node \(i_f\).  These are
the bilinear weights on Cartesian grids and the prescribed sampling weights on
a mesh.  The graph topology remains fixed, but its nonzero interpolation
weights are adapted during transfer fitting.  The deterministic baseline is
\begin{equation}
    \tau_{\ell+1}^{P}(i_f)
    =\sum_{i_c}w^P_{i_fi_c}\tau_\ell(i_c),
    \qquad
    \sum_{i_c}w^P_{i_fi_c}=1.
    \label{eq:resphy_fixed_aggregation}
\end{equation}
At each fine node \(i_f\), the node feature vector is written as
\begin{equation}
    \phi_{i_f}
    =
    \left[
        \tau_{\ell+1}^{P}(i_f),\,
        x_{i_f},\,
        y_{i_f},\,
        \phi_{i_f}^{\rm app},\,
        g_{i_f}^{\rm opt}
    \right]^{T}.
    \label{eq:resphy_node_feature_vector}
\end{equation}
Here \(\phi_{i_f}^{\rm app}\) denotes the application-dependent channels such
as observations, masks, source channels, interpolation weights, or local
coarse-stencil values, and \(g_{i_f}^{\rm opt}\) denotes an optional
loss-diagnostic channel when such a channel is used.
Section~\ref{sec:problem_realizations} gives the Darcy and EIT feature vectors
explicitly.

The shared map \(\operatorname{Corr}_{\omega_\ell}\) is a
two-hidden-layer SiLU multilayer perceptron of width 64, applied pointwise to
\(\phi_{i_f}\).  Concretely,
\[
    \operatorname{Corr}_{\omega_\ell}(\phi)
    =W_2\,\operatorname{SiLU}\bigl(W_1\,\operatorname{SiLU}(W_0\phi+b_0)+b_1\bigr)+b_2,
\]
where \(\omega_\ell=\{W_0,W_1,W_2,b_0,b_1,b_2\}\) and
\(\operatorname{SiLU}(x)=x\,\sigma(x)\) with \(\sigma\) the logistic function.
Its output is a flat vector whose first \(k\) entries are the edge biases and the last entry is the node residual:
\begin{equation}
    (b_{i_fi_c^{(1)}},\ldots,b_{i_fi_c^{(k)}},\operatorname{res}_{i_f})
    =\operatorname{Corr}_{\omega_\ell}(\phi_{i_f}),
    \label{eq:resphy_weight_residual_output}
\end{equation}
where \(k=|\mathcal{N}_c(i_f)|\) and \(\mathcal{N}_c(i_f)=\{i_c^{(1)},\ldots,i_c^{(k)}\}\) is the ordered coarse stencil of fine node \(i_f\).
The output layer is initialized to zero, and a fresh parameter set
\(\omega_\ell\) is fitted at each interface.  The learned aggregation weights
are the interpolation-centered softmax
\begin{equation}
    \widehat w_{i_fi_c}(\omega_\ell)
    =
    \frac{w^P_{i_fi_c}\exp(b_{i_fi_c})}
    {\sum_{j_c}w^P_{i_fj_c}\exp(b_{i_fj_c})},
    \qquad
    \sum_{i_c}\widehat w_{i_fi_c}(\omega_\ell)=1.
    \label{eq:resphy_learned_weights}
\end{equation}
For multichannel representations, separate weight groups may be used for
state and parameter channels; Darcy uses one group for the pressure channels
and one for the raw permeability coordinate, whereas EIT uses one scalar
group.  The learned prolongation and residual update are
\begin{subeqnarray}
    \label{eq:residual_gnn_correction}
    \tau_{\ell+1}^{W}(i_f)
    & = &
    \sum_{i_c}\widehat w_{i_fi_c}(\omega_\ell)\tau_\ell(i_c),\\
    \tau_{\ell+1}(i_f)
    & = &
    \tau_{\ell+1}^{W}(i_f)+\epsilon\tanh(\operatorname{res}_{i_f}).
\end{subeqnarray}
Here \(\tanh\) acts componentwise, \(\epsilon>0\) sets the maximum nodewise
correction amplitude, and application-specific output channels may use
distinct positive amplitudes.
The learned prolongation \(\tau_{\ell+1}^W\) and the incoming representation
\(\tau_{\ell+1}\) may be defined on an auxiliary target
grid.  The fixed incoming map
\(\operatorname{In}_{\ell+1}\) then maps this representation to the next
DNF discretization and converts it into trainable variables, while
\(\operatorname{Asm}_{\ell+1}\) constructs admissible physical fields, restores
prescribed boundary values, and enforces the parameterizations required by
the numerical model.  For Darcy, \(\operatorname{In}_{\ell+1}\) copies the corrected interior pressures
and raw permeability parameter \(\rho\), and
\(\operatorname{Asm}_{\ell+1}\) restores the physical boundary values and
maps \(\rho\) to the positive permeability \(K\).  For the signed EIT contrast,
\(\operatorname{In}_{\ell+1}\) performs the fixed auxiliary-grid-to-mesh
sampling and \(\operatorname{Asm}_{\ell+1}\) is the identity on the mesh
contrast.  At initialization, \(b_{i_fi_c}=0\) and
\(\operatorname{res}_{i_f}=0\), so
\(\widehat w_{i_fi_c}=w^P_{i_fi_c}\) and the complete transfer reduces exactly
to \(\tau_{\ell+1}^P\).  Fixing both outputs at zero defines the deterministic
interpolation ablation; learning only \(b\) gives the learned-weight
ablation, while the complete ResiPhy-GNN learns both weights and residuals.

\begin{figure}[!htbp]
\centering
\includegraphics[width=0.95\linewidth]{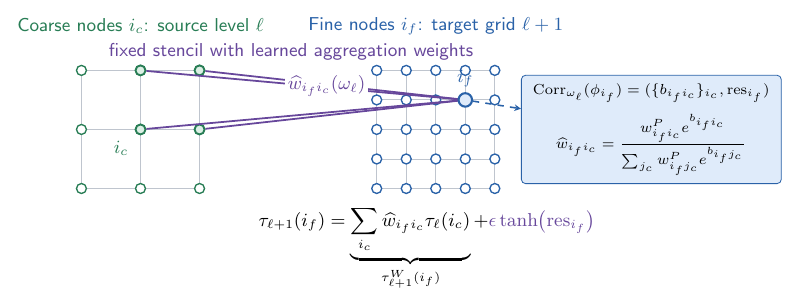}
\caption{Graph-and-node view of ResiPhy-GNN at interface \(\ell\).  The fixed
stencil connects coarse nodes \(i_c\) to fine node \(i_f\).  Its weights are
initialized by \(w^P_{i_fi_c}\), adapted through learned logit biases and
softmax normalization, and used to form \(\tau_{\ell+1}^{W}(i_f)\).  ResiPhy-GNN then
adds a nodewise residual, with \(\omega_\ell\) fitted
by the target-level transfer loss shown in the figure.}
\label{fig:resiphy_gnn_graph_structure}
\end{figure}

Figure~\ref{fig:resiphy_gnn_graph_structure} gives a graph-and-node view of the
coarse-to-refined transfer.

\paragraph{Transfer loss and network fitting.}
Let \(\mathcal T_{\ell+1}(\tau;d)\)
denote the fixed transfer loss on a candidate incoming representation.  It is
evaluated by mapping \(\tau\) through \(\operatorname{In}_{\ell+1}\) and
\(\operatorname{Asm}_{\ell+1}\) to the prescribed target-level numerical
model.  Thus all boundary and admissibility conditions required by the
physical loss blocks are imposed by \(\operatorname{Asm}_{\ell+1}\).  The
transfer loss may also depend directly on \(\tau\), and it need not
equal the reconstruction loss: it may select or reweight data, physics, and
regularization blocks solely to construct an initialization. 
Concretely,
\[
    \mathcal T_{\ell+1}(\tau)
    =\mathcal J^{\rm data}_{\ell+1}(z)
      +w_{\rm pde}\,\mathcal J^{\rm phys}_{\ell+1}(z)
      +w_{\rm reg}\,\mathcal R_{\ell+1}(\tau),
    \qquad
    z=\operatorname{Asm}_{\ell+1}(\operatorname{In}_{\ell+1}(\tau)),
\]
where \(w_{\rm pde},w_{\rm reg}\geq0\) weight the physics and regularization
blocks for transfer (and may differ from the reconstruction
weights in \eqref{eq:mdnf_generic_loss}).
Using the complete reconstruction loss is the special case
\[
    \mathcal T_{\ell+1}(\tau;d)=\mathcal L_{\ell+1}(\operatorname{In}_{\ell+1}(\tau)).
\]
Direct dependence on \(\tau\) permits a transfer-space regularizer such as the EIT
grid-smoothness term.  The data argument \(d\) is suppressed below.

During transfer fitting, \(\omega_\ell\) is optimized with the coarse solution
and transfer loss held constant.  Abstractly, the network parameters solve
\begin{equation}
    \min_{\omega_\ell}\;
    \mathcal T_{\ell+1}
    (\tau_{\ell+1}).
    \label{eq:transfer_training_problem}
\end{equation}
After solving this problem, \(\omega_\ell\) denotes the fitted parameters at
interface \(\ell\), and
\(\tau_{\ell+1}\) in
\eqref{eq:residual_gnn_correction} is the fitted ResiPhy-GNN output.
The loss couples fine nodes through its finite-difference stencil,
finite-element mesh, regularizer, or CEM measurement map, even though the
shared map that predicts edge-logit biases and residuals is applied
pointwise.

\begin{figure}[!htbp]
\centering
\includegraphics[width=0.95\linewidth]{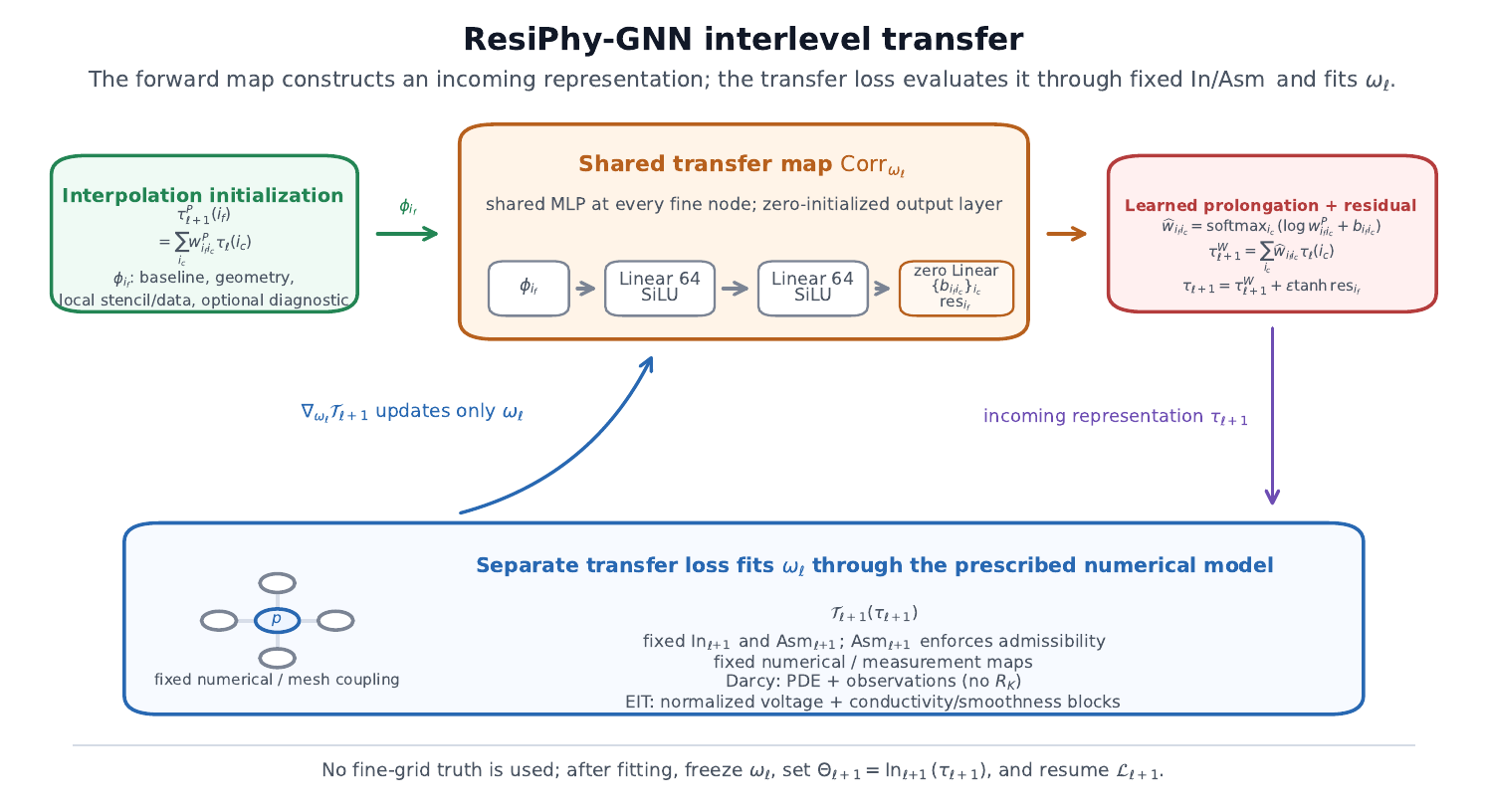}
\caption{ResiPhy-GNN at one coarse-to-fine interface.  The node feature
\(\phi_{i_f}\) enters a shared map that predicts edge-logit biases and a node
residual.  The transfer loss \(\mathcal T_{\ell+1}\), evaluated through the
fixed \(\operatorname{In}_{\ell+1}\), \(\operatorname{Asm}_{\ell+1}\), and
numerical operators, fits \(\omega_\ell\).}
\label{fig:gnn_network_architecture}
\end{figure}

\begin{samepage}
The resulting initialization and fine-level DNF solve are
\begin{subeqnarray}
    \label{eq:fine_solve_after_transfer}
    \Theta_{\ell+1}
    & = &
    \operatorname{In}_{\ell+1}
    \!\left(\tau_{\ell+1}\right),\\
    \Theta_{\ell+1}^{\rm opt}
    & = &
    \argmin_{\Theta_{\ell+1}}
    \mathcal{L}_{\ell+1}(\Theta_{\ell+1}),\\
    z_{\ell+1}^{\rm opt}
    & = & \operatorname{Asm}_{\ell+1}(\Theta_{\ell+1}^{\rm opt}).
\end{subeqnarray}
\end{samepage}

At each interface, \(\omega_\ell\) is reinitialized and fitted for the current
inverse instance.  The Darcy transfer loss uses the refined PDE residual and
observed pressures; the KTC2023 EIT transfer loss uses the linearized voltage
loss and conductivity regularization.

The next section instantiates the abstract operators $\operatorname{Asm}$,
$\operatorname{Out}$, $\operatorname{In}$, the reconstruction
loss $\mathcal L$, and the node feature vector $\phi_{i_f}$ for Darcy
permeability inversion and KTC2023 EIT.
\section{Darcy and EIT Realizations}
\label{sec:problem_realizations}

The governing models and discretizations are defined in
Section~\ref{sec:preliminaries}.  Here we instantiate
Algorithm~\ref{alg:complete_mdnf_gnn} for a full-space Darcy inversion that
jointly optimizes pressure and permeability and for a reduced linearized
complete electrode model (CEM) realization of KTC2023 electrical impedance
tomography (EIT).
Within each realization, we suppress the interface subscript on
\(\omega_\ell\) and write \(\omega\) when only one transfer interface is under
discussion.

\subsection{Darcy full-space realization}
\label{subsec:darcy_realization}

\paragraph{Discrete-field representation and reconstruction loss.}
The Darcy realization uses the full-space formulation: all pressure fields and
their shared permeability remain active variables.  For a level \(\ell\), define
\begin{equation}
\begin{aligned}
    \Theta_\ell
    &=\{U_{1,\ell},\ldots,U_{M,\ell},\rho_\ell\},
    &K_\ell&=K_{\min}+\operatorname{softplus}(\rho_\ell),\\
    \tau_\ell
    &=\operatorname{Out}_\ell(\Theta_\ell^{\rm opt})
      =\Theta_\ell^{\rm opt},
    &\Theta_\ell
    &=\operatorname{In}_\ell(\tau_\ell),\\
    z_\ell
    &=\operatorname{Asm}_\ell(\Theta_\ell)
      =\{U_{1,\ell},\ldots,U_{M,\ell},K_\ell\}.
\end{aligned}
\label{eq:darcy_representation_maps}
\end{equation}
Here \(K_{\min}>0\) and \(\rho_\ell\) is an unconstrained parameter.  The
pointwise map \(\operatorname{softplus}(s):=\log(1+e^s)\) enforces positive
permeability.  The initialization copies transferred interior pressures and
\(\rho_\ell\), while assembly restores the prescribed Dirichlet values and maps
\(\rho_\ell\) to the positive physical permeability \(K_\ell\).  The superscript
\({\rm D}\) is reserved below for the Darcy residual and feature vector.

Let \(\operatorname{Res}^{\rm D}_{m,\ell}(U_{m,\ell},K_\ell)\) denote the fixed conservative residual
from \eqref{eq:prelim_darcy_discrete}.  The Darcy DNF reconstruction loss is
\begin{equation}
    \begin{aligned}
    \mathcal L_\ell(\Theta_\ell)
    &=
    \frac{\lambda_{\rm pde}}{2M}
    \sum_{m=1}^{M}
    \left\|
    \operatorname{Scale}_{m,\ell}\operatorname{Res}^{\rm D}_{m,\ell}(U_{m,\ell},K_\ell)
    \right\|_2^2\\
    &\quad+
    \frac{\lambda_{\rm obs}}{2M}
    \sum_{m=1}^{M}
    \left\|P_{m,\ell}U_{m,\ell}-d_m\right\|_2^2
    +
    \lambda_K R_K(K_\ell).
    \end{aligned}
    \label{eq:dnf_darcy_loss}
\end{equation}
This loss is optimized over \(\Theta_\ell\), with \(K_\ell\) assembled from
\(\rho_\ell\).  Here \(\operatorname{Scale}_{m,\ell}\), \(P_{m,\ell}\), and \(R_K\)
denote the fixed residual normalization, observation operator, and
permeability regularizer; \(\lambda_{\rm pde}\), \(\lambda_{\rm obs}\), and
\(\lambda_K\) are their nonnegative weights.
In the generic decomposition \eqref{eq:mdnf_generic_loss},
\(\mathcal{J}^{\rm data}_\ell\) is the observation term including its
\(\lambda_{\rm obs}\) weight, \(\mathcal{J}^{\rm phys}_\ell\) is the averaged
scaled PDE residual before multiplication by \(\lambda_{\rm pde}\), and
\(\mathcal{R}_\ell=R_K\) with \(\lambda_{\rm reg}=\lambda_K\).
In stacked form, the Darcy residual and scaling operator are
\begin{equation}
    \mathbf F_\ell(z_\ell)
    =
    \begin{bmatrix}
        \operatorname{Res}^{\rm D}_{1,\ell}(U_{1,\ell},K_\ell)\\
        \vdots\\
        \operatorname{Res}^{\rm D}_{M,\ell}(U_{M,\ell},K_\ell)
    \end{bmatrix},
    \qquad
    \operatorname{Scale}_\ell
    =
    M^{-1/2}\operatorname{diag}_{m=1}^{M}
    \operatorname{Scale}_{m,\ell}.
    \label{eq:darcy_stacked_residual}
\end{equation}
Thus the first term of \eqref{eq:dnf_darcy_loss} is
\(\lambda_{\rm pde}\|\operatorname{Scale}_\ell
\mathbf F_\ell(z_\ell)\|_2^2/2\).
\paragraph{ResiPhy-GNN transfer.}
At interface \(\ell\), level \(\ell+1\) is the target discretization; the
symbol \(f_m\) remains reserved for the Darcy source.  Let
\(\mathcal N_c(i_f)=(i_c^{(1)},\ldots,i_c^{(4)})\) denote the ordered
four-point coarse stencil of target node \(i_f\), and let \(w^P_{i_fi_c}\) be
the corresponding prescribed bilinear weights.  A superscript \(P\) denotes
the interpolation baseline:
\begin{subeqnarray}
    \label{eq:darcy_prolongation_baseline}
    U_{m,\ell+1}^{P}(i_f)
    & = & \sum_{q=1}^{4}w^P_{i_fi_c^{(q)}}
        U_{m,\ell}(i_c^{(q)}),\\
    \rho_{\ell+1}^{P}(i_f)
    & = & \sum_{q=1}^{4}w^P_{i_fi_c^{(q)}}
        \rho_{\ell}(i_c^{(q)}),\\
    \tau_{\ell+1}^{P}
    & = & \{U_{1,\ell+1}^P,\ldots,U_{M,\ell+1}^P,
             \rho_{\ell+1}^P\},\\
    K_{\ell+1}^P
    & = & K_{\min}+\operatorname{softplus}(\rho_{\ell+1}^P).
\end{subeqnarray}
The target-node feature vector is
\begin{equation}
    \begin{aligned}
    \phi^{\rm D}_{i_f}
    =\Bigl[&
        x_{i_f},\,
        y_{i_f},\,
        \bigl(w^P_{i_fi_c^{(q)}},\rho_\ell(i_c^{(q)})\bigr)_{q=1}^{4},\,
        K_{\ell+1}^P(i_f),\\[-1mm]
        &\quad
        \bigl(U_{m,\ell+1}^P(i_f),f_{m,\ell+1}(i_f),
        \chi^{\rm obs}_{m,\ell+1}(i_f),d_{m,i_f}\bigr)_{m=1}^{M}
    \Bigr]^{T}
    \in\mathbb R^{4M+11},
    \end{aligned}
    \label{eq:darcy_gnn_features}
\end{equation}
where \((x_{i_f},y_{i_f})\) are target-node coordinates and
\(\chi^{\rm obs}_{m,\ell+1}\), \(f_{m,\ell+1}\), and \(d_{m,i_f}\) are the observation
indicator, source, and observed-pressure channels.  Darcy uses no
loss-gradient feature; physical coupling enters when the transfer map is fitted
through the conservative-residual transfer loss.

For the four-point stencil, the shared transfer map
\(\operatorname{Corr}_\omega:\mathbb R^{4M+11}\to\mathbb R^{M+9}\) has the
ordered output
\begin{subeqnarray}
    \label{eq:darcy_componentwise_gnn_update}
    \operatorname{Corr}_\omega(\phi^{\rm D}_{i_f})
    & = &
    \left[
        \bigl(b^u_{i_fi_c^{(q)}},b^\rho_{i_fi_c^{(q)}}\bigr)_{q=1}^{4},\,
        \bigl(\operatorname{res}^u_{m,i_f}\bigr)_{m=1}^{M},\,
        \operatorname{res}^\rho_{i_f}
    \right]^{T},\\
    \widehat w^s_{i_fi_c}
    & = &
    \frac{w^P_{i_fi_c}\exp(b^s_{i_fi_c})}
    {\sum_{j_c}w^P_{i_fj_c}\exp(b^s_{i_fj_c})},
    \qquad s\in\{u,\rho\},\\
    U_{m,\ell+1}^{W}(i_f)
    & = & \sum_{i_c}\widehat w^u_{i_fi_c}U_{m,\ell}(i_c),\\
    \rho_{\ell+1}^{W}(i_f)
    & = & \sum_{i_c}\widehat w^\rho_{i_fi_c}\rho_{\ell}(i_c),\\
    \widetilde U_{m,\ell+1}(i_f)
    & = &
        U_{m,\ell+1}^{W}(i_f)+\epsilon_u\tanh \operatorname{res}^u_{m,i_f},\\
    \rho_{\ell+1}(i_f)
    & = &
        \rho_{\ell+1}^W(i_f)+\epsilon_\rho\tanh \operatorname{res}^\rho_{i_f}.
\end{subeqnarray}
Here \(\epsilon_u,\epsilon_\rho>0\) set the maximum state and raw
permeability corrections.  The state channels share \(\widehat w^u\), while
the raw permeability uses \(\widehat w^\rho\).  Zero output biases and
residuals recover the bilinear baseline exactly.  The transferred raw
pressures and raw permeability parameter form the incoming representation
\[
    \tau_{\ell+1}
    =\{\widetilde U_{1,\ell+1},\ldots,
       \widetilde U_{M,\ell+1},\rho_{\ell+1}\}.
\]
The admissible physical fields used by the transfer loss are
\begin{equation}
    \{U_{1,\ell+1},\ldots,U_{M,\ell+1},K_{\ell+1}\}
    =\operatorname{Asm}_{\ell+1}\!\left(
        \operatorname{In}_{\ell+1}(\tau_{\ell+1})
    \right),
    \label{eq:darcy_transfer_assembled_fields}
\end{equation}
where \(\operatorname{Asm}_{\ell+1}\) restores the prescribed pressure
boundary values and maps \(\rho_{\ell+1}\) to the positive permeability.

For a candidate \(\tau_{\ell+1}\), the transfer loss is
\begin{equation}
    \mathcal T_{\ell+1}(\tau_{\ell+1})
    =
    w_{\rm pde} \mathcal E_{\rm pde}
    +w_{\rm obs} \mathcal E_{\rm obs},
    \label{eq:warm_start_loss}
\end{equation}
where \(w_{\rm pde},w_{\rm obs}\geq0\).  Let \(\mathcal I_{\ell+1}\) be the interior
target nodes, \(\mathcal O_{m,\ell+1}\) the observed target nodes for source \(m\),
\(N_{\rm int}=M|\mathcal I_{\ell+1}|\), and
\(N_{\rm obs}=\sum_m|\mathcal O_{m,\ell+1}|\).  With the fixed stabilizing floor
\(\varepsilon_{\rm sc}=10^{-12}\), the normalized blocks are
\begin{subeqnarray}
\label{eq:darcy_warm_start_blocks}
 s_{\rm pde}
 & = & \left[
    \frac{1}{N_{\rm int}}
    \sum_{m=1}^M\sum_{i_f\in\mathcal I_{\ell+1}} f_{m,\ell+1}(i_f)^2
    +\varepsilon_{\rm sc}
   \right]^{1/2},\\
 \mathcal E_{\rm pde}
 & = & \frac{1}{N_{\rm int}}
   \sum_{m=1}^M\sum_{i_f\in\mathcal I_{\ell+1}}
   \left[
    \frac{\operatorname{Res}^{\rm D}_{m,\ell+1}(U_{m,\ell+1},K_{\ell+1})(i_f)}{s_{\rm pde}}
   \right]^2,\\
 \mathcal E_{\rm obs}
 & = & \frac{1}{N_{\rm obs}}
   \sum_{m=1}^M\sum_{i_f\in\mathcal O_{m,\ell+1}}
   \bigl(U_{m,\ell+1}(i_f)-d_{m,i_f}\bigr)^2.
\end{subeqnarray}
In the generic form of Section~\ref{subsec:residual_gnn_framework}, this
choice corresponds to
\(\mathcal J^{\rm data}_{\ell+1}=w_{\rm obs}\mathcal E_{\rm obs}\),
\(\mathcal J^{\rm phys}_{\ell+1}=\mathcal E_{\rm pde}\), and \(w_{\rm reg}=0\).
Thus \(\mathcal T_{\ell+1}\), the Darcy realization of the generic transfer
loss, uses only the
target-grid PDE residual and pressure observations; it omits the permeability
regularizer \(R_K\) used by the reconstruction loss.  Network fitting and the
subsequent target-level initialization follow
\eqref{eq:transfer_training_problem} and
\(\Theta_{\ell+1}
=\operatorname{In}_{\ell+1}(\tau_{\ell+1})\), respectively.

\subsection{KTC2023 CEM-EIT reduced-space realization}
\label{subsec:eit_realization}

\paragraph{Discrete-field representation and reconstruction losses.}
Let \(c\), \(a\), and \({\rm mesh}\) denote the coarse coordinate grid, the
auxiliary transfer grid, and the official FEM-node space, respectively.  The
fixed coordinate-to-mesh maps are \(\Pi_c\) and \(\Pi_a\).  Table~\ref{tab:eit_stage_realizations}
summarizes the three representations; the auxiliary stage stores an incoming
transfer representation and is not a DNF level.

\begin{table}[!htbp]
\centering
{\footnotesize
\setlength{\tabcolsep}{3pt}
\renewcommand{\arraystretch}{1.1}
\begin{tabularx}{\linewidth}{@{}l
  >{\centering\arraybackslash}X
  >{\centering\arraybackslash}X
  >{\centering\arraybackslash}X@{}}
\toprule
Stage & Stage representation & FEM field & Loss function \\
\midrule
Coarse DNF &
\(\Theta_c=\tau_c=\vartheta_c\) &
\(z_c=\Pi_c\vartheta_c\) &
\(\mathcal L_c\) \\
Auxiliary transfer &
\(\tau_a\) &
\(\Pi_a\tau_a\) &
\(\mathcal T_a\) \\
Mesh DNF &
\(\Theta_{\rm mesh}=\delta\sigma_{\rm mesh}\) &
\(z_{\rm mesh}=\delta\sigma_{\rm mesh}\) &
\(\mathcal L_{\rm mesh}\) \\
\bottomrule
\end{tabularx}
}
\caption{Representations in the EIT realization.  The auxiliary transfer is
an interface representation between the two DNF levels.}
\label{tab:eit_stage_realizations}
\end{table}

Here \(\operatorname{Out}_c\) is the identity on \(\vartheta_c\),
\(\operatorname{Asm}_c=\Pi_c\),
\(\operatorname{In}_{\rm mesh}=\Pi_a\), and
\(\operatorname{Asm}_{\rm mesh}=I_{\rm mesh}\).  Applying the fixed
aggregation \eqref{eq:resphy_fixed_aggregation} from the coarse grid to the
auxiliary grid, with fixed edge set \(\mathsf E_c^a\), gives
\begin{equation}
    \tau_a^{P}(i_f)
    =\sum_{(i_c,i_f)\in\mathsf E_c^a}w^P_{i_f i_c}\tau_c(i_c),
    \qquad
    \sum_{(i_c,i_f)\in\mathsf E_c^a}w^P_{i_f i_c}=1.
    \label{eq:eit_auxiliary_baseline}
\end{equation}

For any stage coordinate \(q\) mapped to the FEM mesh by a fixed linear map
\(\Pi\), define the voltage, conductivity, and grid-smoothness blocks
\begin{subeqnarray}
    \label{eq:eit_loss_blocks}
    \mathcal{E}_{V}(q;\Pi)
    & = & \frac{1}{2N}
    (J\Pi q-\Delta V)^T W(J\Pi q-\Delta V),\\
    \mathcal{E}_{\sigma}(q;\Pi)
    & = & \frac{1}{2N_{\rm node}}\|L_\sigma\Pi q\|_2^2,\\
    \mathcal{E}_{\rm sm}(q)
    & = & R_{\rm grid}(q).
\end{subeqnarray}
Here \(W\) is the full measurement precision matrix, \(N\) is the number of
active voltage measurements, \(N_{\rm node}\) is the number of FEM nodes,
\(L_\sigma\) is the fixed conductivity regularization operator, and
\(R_{\rm grid}\) is the regular-grid smoothness functional.  With
\(I_{\rm mesh}\) denoting the FEM identity, the three stage losses are
\begin{subeqnarray}
    \label{eq:eit_linear_objective}
    \mathcal L_c(\vartheta_c)
    & = & \mathcal{E}_V(\vartheta_c;\Pi_c)
      +\mathcal{E}_\sigma(\vartheta_c;\Pi_c)
      +\lambda_{\rm sm}\mathcal{E}_{\rm sm}(\vartheta_c),\\
    \mathcal T_a(\tau_a)
    & = & \frac{\mathcal{E}_V(\tau_a;\Pi_a)}
                    {s_V^{\rm tr}}
      +\lambda_\sigma
       \frac{\mathcal{E}_\sigma(\tau_a;\Pi_a)
             +\lambda_{\rm sm}\mathcal{E}_{\rm sm}(\tau_a)}
            {s_\sigma^{\rm tr}},\\
    \mathcal L_{\rm mesh}(\delta\sigma_{\rm mesh})
    & = & \frac{\mathcal{E}_V(\delta\sigma_{\rm mesh};I_{\rm mesh})}
                    {s_V^{\rm mesh}}
      +\lambda_\sigma
       \frac{\mathcal{E}_\sigma(\delta\sigma_{\rm mesh};I_{\rm mesh})}
            {s_\sigma^{\rm mesh}}.
\end{subeqnarray}
The nonnegative weights are stage specific: \(\lambda_{\rm sm}\) scales grid
smoothness, and \(\lambda_\sigma\) scales the normalized conductivity block
and its transfer-grid smoothness term.  The three functions in
\eqref{eq:eit_linear_objective} are the coarse reconstruction, interface
transfer, and mesh reconstruction losses, respectively.  In the generic
decomposition \eqref{eq:mdnf_generic_loss},
\(\mathcal{J}^{\rm data}_h=\mathcal{E}_V\),
\(\mathcal{J}^{\rm phys}_h=0\) (the state is eliminated by the linearized
CEM), and \(\mathcal{R}_h\) aggregates
\(\mathcal{E}_\sigma+\lambda_{\rm sm}\mathcal{E}_{\rm sm}\), up to the
stage-dependent normalization factors.

The transfer scales \(s_V^{\rm tr}\) and \(s_\sigma^{\rm tr}\) are the positive
values of the voltage and combined conductivity--smoothness blocks at the
uncorrected auxiliary-grid baseline.  The mesh scales
\(s_V^{\rm mesh}\) and \(s_\sigma^{\rm mesh}\) are the corresponding block
values after this baseline is mapped to the FEM mesh.  Consequently, each
normalized block starts at unit scale in the stage where it is used.

The conductivity weight is selected from the measured voltage energy:
\begin{subeqnarray}
 \label{eq:eit_energy_lambda}
 S_{\rm data}
 & = & \frac{1}{2N}\Delta V^T W\Delta V,\\
 \lambda_\sigma
 & = & \operatorname{clip}
 \left(c_S S_{\rm data},\lambda_{\min},\lambda_{\max}\right).
\end{subeqnarray}
Here \(\operatorname{clip}(x,a,b):=\min\{\max\{x,a\},b\}\), and
\(c_S>0\), \(\lambda_{\min}\), and \(\lambda_{\max}\) are fixed constants
reported in Section~\ref{subsec:experiment_setup}.  No reconstruction or
segmentation labels enter their calibration or the evaluation of
\eqref{eq:eit_energy_lambda}.

The optimized variable is an unconstrained signed contrast, not an absolute
conductivity, so no positivity or box projection is applied.  Because \(J\) is
assembled at the positive background \(\sigma_0\), negative and positive
contrasts represent resistive and conductive inclusions; the linearized model
does not require \(\sigma_0+\delta\sigma\) to remain positive.  Other KTC2023
approaches use spatial or level-set priors, hierarchical Bayesian inversion,
and multibang regularization
\cite{alghamdi2024spatial,pragliola2024hierarchical,brandt2024multibang};
here the official linearized CEM model is retained to isolate the multilevel
discrete-field path.

\paragraph{ResiPhy-GNN transfer.}
At auxiliary node \(i_f\), let \(N_a\) denote the number of auxiliary-grid
nodes.  The transfer-loss gradient diagnostic and its fixed root-mean-square
scales are
\begin{subeqnarray}
    \label{eq:eit_feature_scales}
    g_{i_f}^{\rm EIT}
    & = & \left.
        \nabla_{\tau_a(i_f)}
        \mathcal T_a(\tau_a)
        \right|_{\tau_a=\tau_a^{P}},\\
    s_\tau
    & = & \max\!\left\{
        \left[\frac{1}{N_a}\sum_{i_f=1}^{N_a}
        \bigl(\tau_a^{P}(i_f)\bigr)^2
        \right]^{1/2},1\right\},\\
    s_g
    & = & \max\!\left\{
        \left[\frac{1}{N_a}\sum_{i_f=1}^{N_a}(g_{i_f}^{\rm EIT})^2
        \right]^{1/2},10^{-12}\right\}.
\end{subeqnarray}
The EIT feature vector is
\begin{equation}
    \phi^{\rm EIT}_{i_f} =
    \left(
        x_{i_f},y_{i_f},
        \frac{\tau_a^{P}(i_f)}{s_\tau},
        \frac{g_{i_f}^{\rm EIT}}{s_g}
    \right),
    \label{eq:eit_gnn_features}
\end{equation}
where \((x_{i_f},y_{i_f})\in[-1,1]^2\) are normalized auxiliary-grid coordinates.
The interpolation baseline is a node feature, whereas \(g_{i_f}^{\rm EIT}\) is the
optional loss diagnostic defined in \eqref{eq:eit_feature_scales}, not a
message-passing operation.

The five outputs of the shared transfer map define four local weight biases
and one node residual.  The learned aggregation and incoming representation are
\begin{subeqnarray}
    \label{eq:eit_gnn_update}
    (\{b^\sigma_{i_fi_c}\}_{i_c},\operatorname{res}_{i_f}^\sigma)
    & = & \operatorname{Corr}_\omega(\phi_{i_f}^{\rm EIT}),\\
    \widehat w^\sigma_{i_fi_c}
    & = &
    \frac{w^P_{i_fi_c}\exp(b^\sigma_{i_fi_c})}
    {\sum_{j_c}w^P_{i_fj_c}\exp(b^\sigma_{i_fj_c})},\\
    \tau_a^W(i_f)
    & = & \sum_{i_c}\widehat w^\sigma_{i_fi_c}\tau_c(i_c),\\
    \tau_a(i_f)
    & = & \tau_a^{W}(i_f)
      +\epsilon_\sigma\tanh \operatorname{res}^\sigma_{i_f} ,
\end{subeqnarray}
where \(\epsilon_\sigma>0\) sets the maximum contrast correction.  Zero
weight biases and residual recover \(\tau_a^P\) exactly.  The parameters
\(\omega\) are fitted by \eqref{eq:transfer_training_problem} with
\(\mathcal T_{\ell+1}=\mathcal T_a\).  The maps in
Table~\ref{tab:eit_stage_realizations} then initialize the mesh DNF, while
\(J\), \(W\), the active voltage mask, and \(L_\sigma\) remain fixed.
\section{Numerical Experiments}
\label{sec:experiments}

We use a basic Darcy inverse problem to show that repeated single-level DNF
optimization along a coarse-to-fine path attains higher final accuracy after
the same number of final-grid updates.  We then isolate the inter-level transfer
on LUNDIsim and compare deterministic interpolation with ResiPhy-GNN correction,
followed by a complete two-dimensional four-level ResiPhy-MDNF realization on a
complex Darcy field.  The section closes
with KTC2023 EIT,
where the official complete electrode model linearized reconstruction provides
a strong discrete inverse baseline.

In the comparisons below, the final grid is the discretization associated with
the final coefficient or contrast reconstruction space \(V_J\) defined in
Section~\ref{subsec:mdnf_backbone}.

\subsection{Experimental setup}
\label{subsec:experiment_setup}

The Darcy experiments use a conservative finite-difference residual.  The
unknowns are the state field \(U\) and the permeability coefficient \(K\), and
both are represented directly as trainable discrete fields.  The basic Darcy
case uses a manufactured coefficient to isolate the effect of the multilevel
path under a controlled coefficient field.  LUNDIsim replaces the
manufactured coefficient by reservoir-like fields
\cite{duval2025lundisim}.  The solve on each level optimizes the stated DNF
reconstruction loss.  Each Darcy target uses a newly initialized corrector
fitted to the refined PDE residual and observed pressure mismatch;
\(U^\star\) and \(K^\star\) are reserved for evaluation.  We use
\(w_{\rm pde}=1\), \(w_{\rm obs}=1000\), and
\(\epsilon_u=\epsilon_\rho=0.1\); the correctors are fitted for
3000 Adam steps.

The KTC2023 EIT experiments use the official complete electrode model
finite-element Jacobian and the official active voltage mask for category 1.
For each of the four phantoms, we optimize the linearized difference loss
defined in \eqref{eq:eit_linear_objective}.  We report the relative voltage
residual \(\mathrm{relV}=\|J\delta\sigma-\Delta V\|_2/\|\Delta V\|_2\)
and the official two-threshold Otsu mIoU over background, resistive, and
conductive classes.  Since Otsu thresholding is a postprocessing step, we also
use KTC2023 to examine whether improvements in the continuous reconstruction
are reflected in the thresholded three-class mIoU.
For every phantom, the transfer corrector is reinitialized and trained only
through the fixed linearized loss.  No conductivity truth is
used for that training, and no learned weights are shared between phantoms or
with the Darcy experiments.
The measurement-energy rule \eqref{eq:eit_energy_lambda} uses
\(c_S=8.7228432763\times10^{-6}\), \(\lambda_{\min}=0.01\), and
\(\lambda_{\max}=0.3\).  The constant \(c_S\) is calibrated once from the
Data~1 voltage energy so that the upper cap is attained.  For Data~1--4, the
resulting weights are
\[
    \bigl(
      \lambda_\sigma^{(1)},\lambda_\sigma^{(2)},
      \lambda_\sigma^{(3)},\lambda_\sigma^{(4)}
    \bigr)
    =
    (0.3000,\,0.1464,\,0.1023,\,0.04195).
\]

The coordinate-network PINN baseline uses all 76 current-injection patterns,
a four-layer width-128 network, 10000 Adam updates, and GPU training.  Its
conductivity and state networks are fitted from the PDE, CEM boundary, current,
and measured-voltage losses without the official FEM Jacobian.  The same
three-class IoU averaging convention is used for evaluation.  On Data~3 and
Data~4, however, the recovered field does not admit a valid two-threshold Otsu
split.  Those entries, and consequently the aggregate PINN mean, are reported
as N/A rather than mixing a binary fallback with the required three-class
metric.

The direct single-level baseline starts from zero on one \(64^2\) reconstruction
grid and performs 10000 Adam updates with the full measurement matrix \(W\)
and the conductivity-regularization setting used by the official
reconstruction.  It has no coarse solve, inter-level transfer, ResiPhy-GNN, or
mesh-level refinement.  In the notation of Section~\ref{subsec:eit_realization},
the complete ResiPhy-MDNF pipeline uses a \(32^2\) coarse grid \(c\), a
\(64^2\) auxiliary transfer grid \(a\), and a 1602-node FEM-mesh DNF level.  It
first performs 1200 Adam updates
on the coarse grid.  One ResiPhy-GNN transfer, fitted for eight L-BFGS
iterations using step length \(\eta_k=0.3\), then constructs the incoming
representation on grid \(a\); the fixed incoming map
\(\operatorname{In}_{\mathrm{mesh}}\) samples this corrected field to
the FEM mesh.  There is no intervening
\(64^2\) DNF solve.  Finally, the mesh DNF performs 60 outer L-BFGS updates with
full \(W\).  Thus the single-level row tests direct one-level DNF
optimization.  The causal transfer ablation is instead provided by
the matched LUNDIsim experiment in Section~\ref{subsec:transfer_ablation}; the
KTC2023 experiment compares complete reconstruction pipelines across numerical,
coordinate-network, single-level discrete-field, and ResiPhy-MDNF approaches.

\subsection{Basic Darcy: higher final accuracy under equal final-grid steps}
\label{subsec:basic_darcy_c2f}

Table~\ref{tab:basic_darcy_c2f} reports the controlled basic Darcy hierarchy
and its matched direct \(128^2\) comparison.  The two paths use the same final grid,
6000 Adam steps at \(128^2\), 16 sources, an observation fraction of 0.35,
learning rate \(5\times10^{-4}\), loss weights, and coefficient-freeze
schedule.  The direct solve reaches coefficient and state errors of
\(5.4519\times10^{-3}\) and \(4.7396\times10^{-4}\), respectively.  The
\(64^2\!\to128^2\) path instead reaches \(8.0597\times10^{-4}\) and
\(4.5566\times10^{-5}\), which are \(6.76\times\) and \(10.4\times\) lower.
The later \(256^2\) and \(512^2\) levels mainly preserve the recovered
coefficient while maintaining the state error at the same order.

For a multilevel run completed through level \(\ell\), we define the
cumulative grid-work proxy
\begin{equation}
    C_{{\rm grid},\ell}
    =
    \sum_{k=0}^{\ell}s_k n_k^2,
    \qquad
    \widehat C_{{\rm grid},\ell}
    =
    \frac{C_{{\rm grid},\ell}}{6000\cdot128^2},
    \label{eq:darcy_grid_work}
\end{equation}
where \(\ell\) is the finest completed level, \(s_k=6000\) is the number of
Adam steps at level \(k\), and \(n_k\) is the number of grid cells per
coordinate direction, so \(n_k^2\) is the two-dimensional grid size.  The sum
includes all completed levels \(k=0,\ldots,\ell\).  The denominator is the cost
of the direct single-level \(128^2\) run.  The number of sources and trainable
field channels is common to both paths and therefore cancels in the
normalization.  This proxy measures grid-point operations, not wall-clock time.
For a uniform accuracy comparison, we also define the direct-to-level error
ratios
\begin{equation}
    \mathrm{ER}_{K,\ell}
    =\frac{E_K^{\rm direct}}{E_{K,\ell}},
    \qquad
    \mathrm{ER}_{U,\ell}
    =\frac{E_U^{\rm direct}}{E_{U,\ell}}.
    \label{eq:darcy_error_ratio}
\end{equation}
Here \(\mathrm{ER}\) denotes the error ratio, and \(E_{K,\ell}\) and
\(E_{U,\ell}\) are the final relative \(L^2\) errors of the coefficient and
state fields after the DNF solve on \(V_\ell\).  The superscript ``direct''
denotes the matched direct \(128^2\) baseline. Thus, a ratio greater than one
denotes an error reduction relative to the direct \(128^2\) baseline, whereas
a ratio below one denotes a larger error.

\begin{table}[!htbp]
\centering
{\footnotesize
\setlength{\tabcolsep}{3pt}
\renewcommand{\arraystretch}{1.1}
\begin{tabularx}{\linewidth}{@{}>{\raggedright\arraybackslash}Xccccc@{}}
\toprule
Path / level & \(\widehat C_{{\rm grid},\ell}\) &
Initial \(E_K\) & Final \(E_K\) &
Final \(E_U\) &
\(\mathrm{ER}_{K,\ell}/\mathrm{ER}_{U,\ell}\) \\
\midrule
Direct single-level \(128^2\) & \(1.00\) & \(4.4744\times10^{-1}\) & \(5.4519\times10^{-3}\) &
\(4.7396\times10^{-4}\) & \(1.00/1.00\) \\
Multilevel \(64^2\) & \(0.25\) & \(4.4707\times10^{-1}\) & \(8.7982\times10^{-3}\) &
\(1.3499\times10^{-3}\) & \(0.62/0.35\) \\
Multilevel \(64^2\!\to128^2\) & \(1.25\) & \(6.8516\times10^{-3}\) & \(8.0597\times10^{-4}\) &
\(4.5566\times10^{-5}\) & \(6.76/10.4\) \\
Multilevel \(\to256^2\) & \(5.25\) & \(6.2484\times10^{-4}\) & \(6.2484\times10^{-4}\) &
\(4.1190\times10^{-5}\) & \(8.73/11.5\) \\
Multilevel \(\to512^2\) & \(21.25\) & \(5.4385\times10^{-4}\) & \(5.4385\times10^{-4}\) &
\(4.1763\times10^{-5}\) & \(10.0/11.3\) \\
\bottomrule
\end{tabularx}
}
\caption{Basic Darcy direct and coarse-to-fine comparison.  The work proxy
\(\widehat C_{{\rm grid},\ell}\) is normalized by the direct \(128^2\) run, and
\(\mathrm{ER}>1\) denotes a lower error than that baseline.}
\label{tab:basic_darcy_c2f}
\end{table}

The matched \(128^2\) comparison changes the optimization path, not the final
discretization or loss function.  The additional 6000-step \(64^2\) solve raises
the cumulative work proxy from \(1.00\) to \(1.25\).  For this 25\% increase,
the coefficient and state errors are reduced by factors of \(6.76\) and
\(10.4\), respectively.  This is a quantified total-work comparison under
\eqref{eq:darcy_grid_work}.

\subsection{Transfer ablation: learned transfer improves refined reconstruction}
\label{subsec:transfer_ablation}

After the multilevel backbone is fixed, we isolate the transfer from a coarse
solve to a refined solve.  Table~\ref{tab:transfer_ablation} compares three
\(64^2\to128^2\) LUNDIsim initializations in which the fine-grid coefficient
remains trainable.  The learned transfers are fitted to this target instance
using the physics-and-observation warm-start loss in
\eqref{eq:warm_start_loss}; the fine-grid reference fields are withheld until
evaluation.  After the same refined optimization, the proposed
learned-weights-plus-residual transfer improves both final errors over simple
interpolation: the final state and coefficient errors are approximately
\(9\%\) and \(8\%\) lower, respectively.  The learned-weights-only variant is
included as an ablation to isolate the contribution of the residual update.

\begin{table}[!htbp]
\centering
{\footnotesize
\setlength{\tabcolsep}{3pt}
\renewcommand{\arraystretch}{1.1}
\begin{tabularx}{\linewidth}{@{}>{\raggedright\arraybackslash}Xcc@{}}
\toprule
Fine-grid initialization &
\shortstack{Final \(U\) error\\\((\times10^{-4})\)} &
\shortstack{Final \(K\) error\\\((\times10^{-2})\)} \\
\midrule
Simple interpolation & \(1.8092\) & \(1.5031\) \\
Learned weights only & \(1.6958\) & \(1.5550\) \\
ResiPhy-GNN (learned weights \(+\) residual) & \(1.6456\) & \(1.3787\) \\
\bottomrule
\end{tabularx}
}
\caption{Transfer ablation on the two-dimensional LUNDIsim
\(64^2\to128^2\) Darcy task.  The refined objective and optimization budget
are identical across rows; only the initialization changes.}
\label{tab:transfer_ablation}
\end{table}

\begin{figure}[!htbp]
\centering
\includegraphics[width=0.95\linewidth]{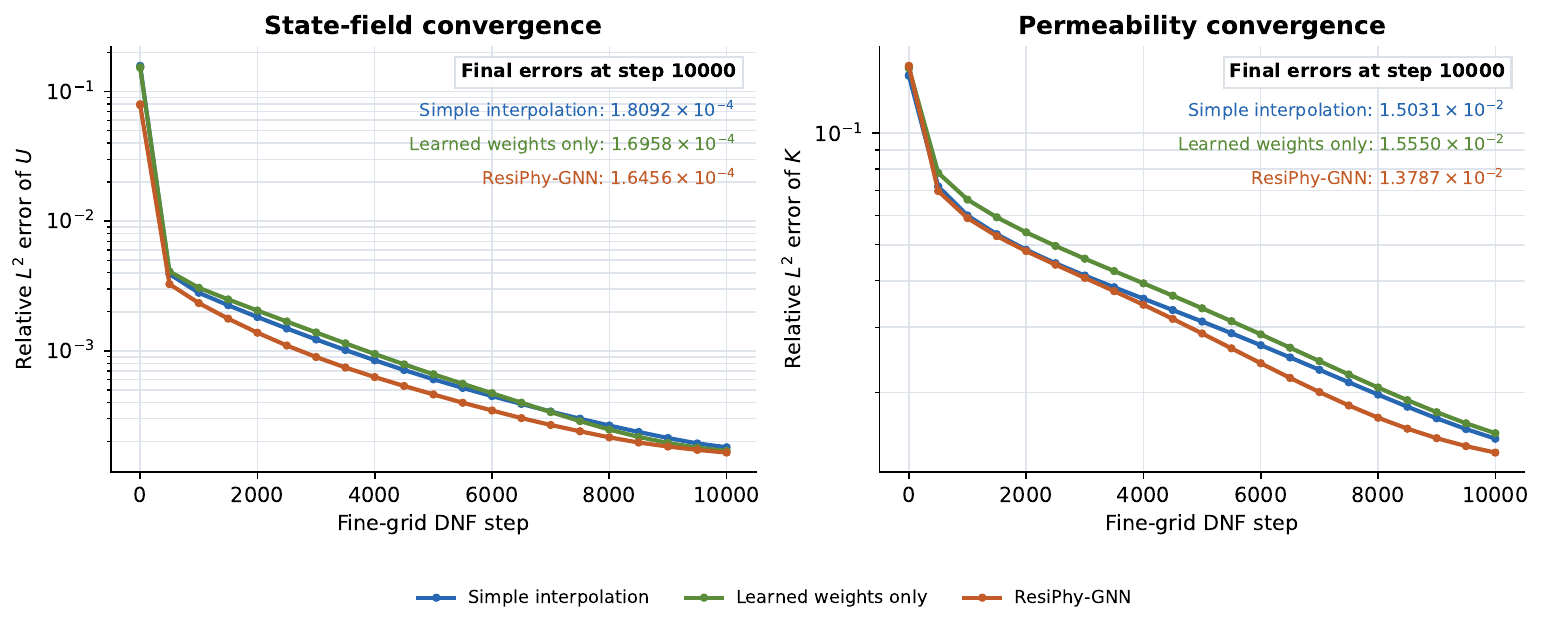}
\caption{Aligned fine-grid convergence for the three transfer initializations
in Table~\ref{tab:transfer_ablation} on the two-dimensional LUNDIsim problem.
All three receive the same 10000 target-grid DNF updates; step zero denotes
the transferred initialization.  The preceding 3000 transfer-training steps
for the two learned initializations are not included on this horizontal axis.}
\label{fig:gnn_aligned_fine_convergence}
\end{figure}

Figure~\ref{fig:gnn_aligned_fine_convergence} compares all three rows of
Table~\ref{tab:transfer_ablation}, aligning their trajectories by the
number of target-grid DNF updates after their respective initializations have
been constructed.  The learned-weights-plus-residual initialization starts
with a substantially lower state error but a slightly higher coefficient
error.  After the first recorded 500 fine-grid steps, both of its error curves
are below their simple-interpolation counterparts and retain this advantage
through 10000 steps.  Learning the weights alone ultimately lowers the state
error but not the coefficient error.  At the common endpoint,
Table~\ref{tab:transfer_ablation} gives state
errors of \(1.6456\times10^{-4}\) versus \(1.8092\times10^{-4}\) and
coefficient errors of \(1.3787\times10^{-2}\) versus
\(1.5031\times10^{-2}\).  Thus the learned transfer yields a better refined
optimization path under the same number of fine-grid updates.  The 3000
transfer-training steps occur before the axis in
Figure~\ref{fig:gnn_aligned_fine_convergence}.

\subsection{LUNDIsim extension: complete multilevel inversion of a complex Darcy field}
\label{subsec:lundisim_extensions}

The preceding transfer ablation establishes the matched two-dimensional
comparison on a reservoir-like LUNDIsim coefficient, so we do not repeat those
numbers here.  We next record a longer multilevel realization that extends the
scope of the Darcy result.  This is a complete
\(16^2\to32^2\to64^2\to128^2\) ResiPhy-MDNF run: after the \(16^2\) DNF solve,
a newly initialized ResiPhy-GNN is fitted at each of the three interlevel
interfaces using \(\mathcal T_{\ell+1}\) in
\eqref{eq:warm_start_loss}.  Each corrector receives 3000 Adam updates, and
each of the four DNF levels receives 10000 Adam updates.

Table~\ref{tab:lundisim_complete_transfers} records the transfer-stage PDE
diagnostic before and after each fitted correction.  ResiPhy-GNN reduces the
normalized target-grid PDE loss at all three interfaces, with the largest
factor occurring at the final \(64^2\to128^2\) transition.

\begin{table}[!htbp]
\centering
{\footnotesize
\setlength{\tabcolsep}{3pt}
\renewcommand{\arraystretch}{1.1}
\begin{tabularx}{\linewidth}{@{}c
  >{\centering\arraybackslash}X
  >{\centering\arraybackslash}X
  c@{}}
\toprule
Transition & Bilinear baseline \(\mathcal E_{\rm pde}\) & ResiPhy-GNN \(\mathcal E_{\rm pde}\) & Reduction \\
\midrule
\(16^2\to32^2\)  & \(9.0845\times10^{-1}\) & \(7.0117\times10^{-1}\) & \(1.30\times\) \\
\(32^2\to64^2\)  & \(6.0200\times10^{-1}\) & \(4.0851\times10^{-1}\) & \(1.47\times\) \\
\(64^2\to128^2\) & \(3.7002\times10^{-1}\) & \(1.4700\times10^{-1}\) & \(2.52\times\) \\
\bottomrule
\end{tabularx}
}
\caption{The three fitted ResiPhy-GNN transfers in the complete LUNDIsim
ResiPhy-MDNF realization.  Values are evaluated on the incoming target grid
before the subsequent DNF solve.}
\label{tab:lundisim_complete_transfers}
\end{table}

For this experiment, the relative residual RMS is normalized by the root mean
square of the discrete source:
\begin{equation}
    E_R=
    \frac{
        \sqrt{\operatorname{mean}_{m,i,j}
        \left(\operatorname{Res}^{\rm D}_{m,h}(U_{m,h},K_h)[i,j]\right)^2}
    }{
        \sqrt{\operatorname{mean}_{m,i,j}
        \left(f_{m,h}[i,j]\right)^2}
    },
    \label{eq:lundisim_relative_residual_rms}
\end{equation}
where both means are taken over all sources and interior grid nodes.  After the
final \(128^2\) DNF solve, the complete pipeline reaches a coefficient error of
\(1.3906\times10^{-2}\), a state error of \(1.6766\times10^{-4}\), and
\(E_R=1.4913\times10^{-3}\).
\begin{figure}[!htbp]
\centering
\begin{subfigure}{0.31\linewidth}
    \centering
    \includegraphics[width=\linewidth]{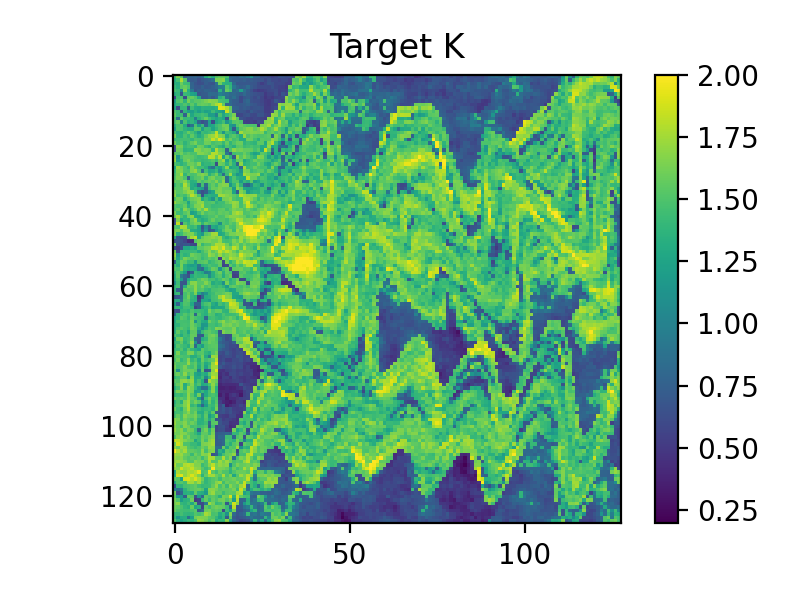}
    \caption{Target \(K\)}
\end{subfigure}
\hfill
\begin{subfigure}{0.31\linewidth}
    \centering
    \includegraphics[width=\linewidth]{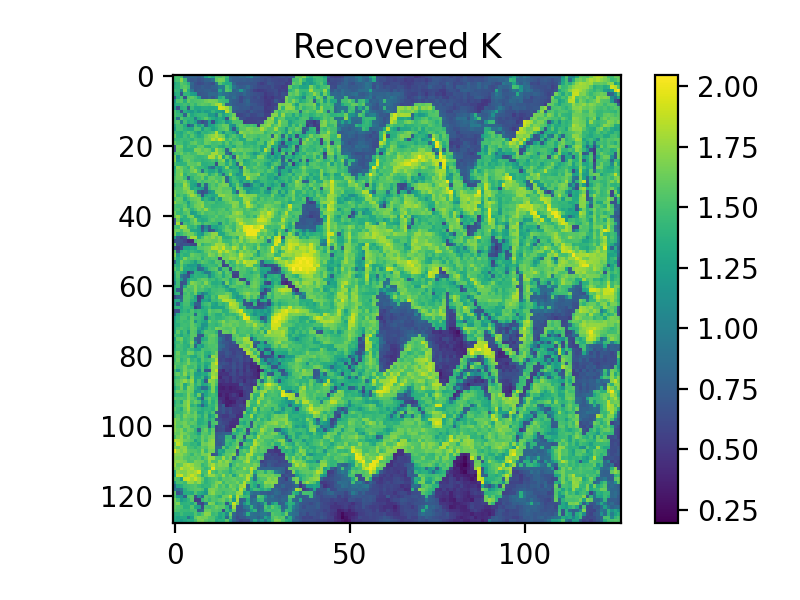}
    \caption{Recovered \(K\)}
\end{subfigure}
\hfill
\begin{subfigure}{0.31\linewidth}
    \centering
    \includegraphics[width=\linewidth]{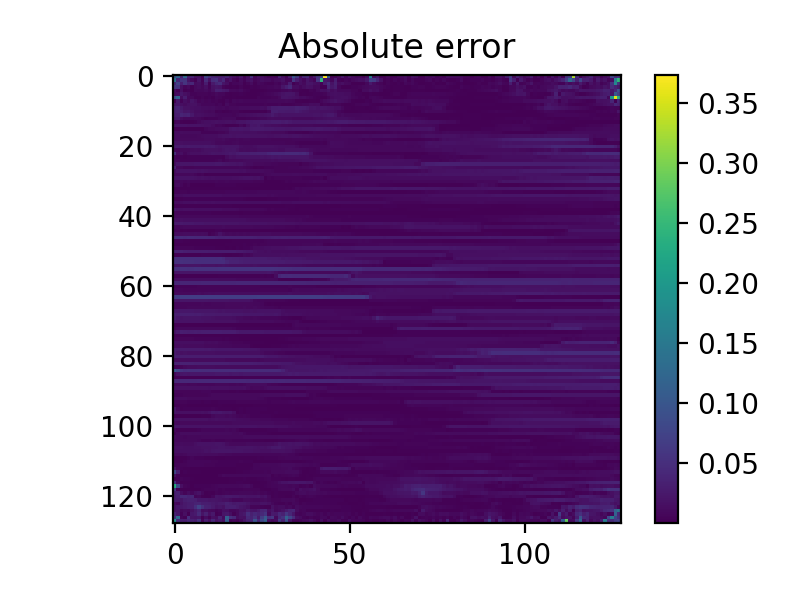}
    \caption{Absolute error}
\end{subfigure}
\caption{Representative two-dimensional LUNDIsim recovery obtained by the
complete \(16^2\to32^2\to64^2\to128^2\) ResiPhy-MDNF pipeline, with a newly
fitted ResiPhy-GNN at each of the three interlevel interfaces.}
\label{fig:lundisim_k_visual}
\end{figure}

\subsection{KTC2023 EIT: cross-physics framework realization}
\label{subsec:ktc2023_pipeline}

KTC2023 changes both the inverse problem and the discretization: the forward
model is the complete electrode model, and the official baseline is a
linearized finite-element reconstruction.  We instantiate the complete
ResiPhy-MDNF workflow with this new physics, discretization, and observation
geometry; all corrector parameters are reinitialized rather than transferred
from Darcy.  Table~\ref{tab:ktc2023_pipeline} compares the official
linearized CEM reconstruction, the evaluated coordinate-network PINN, direct
single-level DNF, and the complete ResiPhy-MDNF pipeline containing ResiPhy-GNN
transfer and mesh-level refinement.  The complete pipeline uses the full
\(W\) matrix and the measurement-energy rule in
\eqref{eq:eit_energy_lambda}; no phantom-specific manual selection of
\(\lambda_\sigma\) is performed.  The transfer corrector uses eight L-BFGS
iterations with step length \(\eta_k=0.3\) for every phantom.  The complete pipeline
reaches \(0.623\) mean mIoU and \(0.1280\) mean relative voltage residual.  The direct single-level
DNF remains the official-setting baseline and is not part of this adaptive
multilevel weighting strategy.

\begin{table}[!htbp]
\centering
{\footnotesize
\setlength{\tabcolsep}{3pt}
\renewcommand{\arraystretch}{1.1}
\begin{tabularx}{\linewidth}{@{}>{\raggedright\arraybackslash}Xcccccc@{}}
\toprule
Method & Mean mIoU & Mean \(\mathrm{relV}\) & Data 1 & Data 2 & Data 3 & Data 4 \\
\midrule
Official linearized CEM baseline & \(0.603\) & -- &
\(0.745\) & \(0.764\) & \(0.506\) & \(0.395\) \\
Coordinate-network PINN, 76 patterns & N/A & -- &
\(0.201\) & \(0.159\) & N/A & N/A \\
Single-level DNF & \(0.533\) & \(0.1367\) &
\(0.725\) & \(0.666\) & \(0.417\) & \(0.325\) \\
Complete ResiPhy-MDNF, mesh-level DNF refinement, full \(W\) & \(0.623\) & \(0.1280\) &
\(0.792\) & \(0.780\) & \(0.509\) & \(0.411\) \\
\bottomrule
\end{tabularx}
}
\caption{KTC2023 EIT comparison of the official reconstruction, coordinate
PINN, direct DNF, and complete ResiPhy-MDNF.  All reported mIoU values use the
same three-class convention; N/A denotes the absence of a valid two-threshold
Otsu split.}
\label{tab:ktc2023_pipeline}
\end{table}

\begin{figure}[!htbp]
\centering
\includegraphics[width=0.95\linewidth]{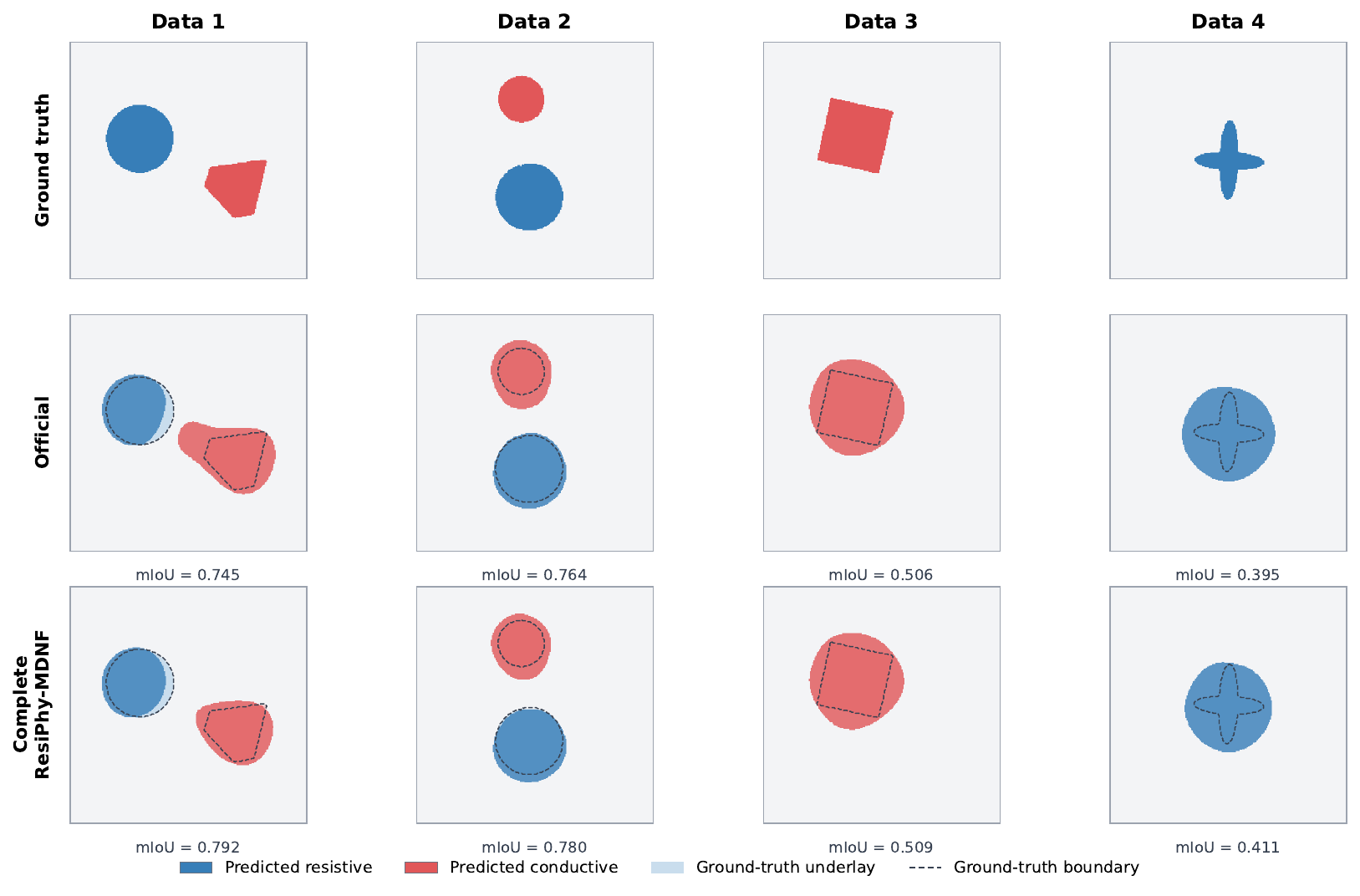}
\caption{Final KTC2023 three-class segmentations.  Columns are Data~1--4 and
rows are the ground truth, official reconstruction, and complete
ResiPhy-MDNF.  Pale underlays and dashed contours indicate ground-truth
supports; mIoU uses the protocol of Table~\ref{tab:ktc2023_pipeline}.}
\label{fig:ktc2023_segmentation_comparison}
\end{figure}

Figure~\ref{fig:ktc2023_segmentation_comparison} shows the final Otsu
segmentations underlying the official and complete-pipeline entries in
Table~\ref{tab:ktc2023_pipeline}.  The complete ResiPhy-MDNF reconstruction
retains the correct anomaly classes in all four phantoms and improves the
three-class mIoU in every case.  The gains on Data~3 and Data~4 are modest,
consistent with the sensitivity of their sharp square and cross supports to
thresholding and boundary smoothing.

\begin{figure}[!htbp]
\centering
\includegraphics[width=0.95\linewidth]{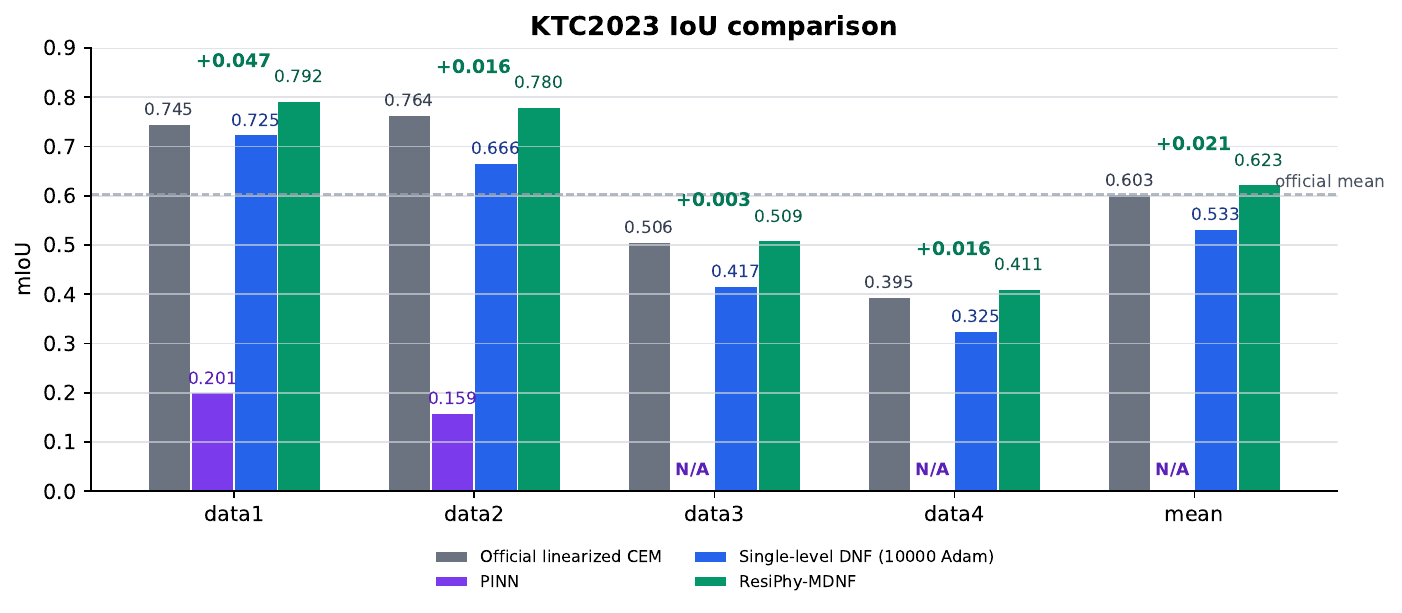}
\caption{KTC2023 mIoU comparison of the official CEM baseline, coordinate
PINN, direct DNF, and complete ResiPhy-MDNF.  Signed annotations give the
complete-pipeline difference from the official baseline; unavailable PINN
segmentations are marked N/A.}
\label{fig:ktc2023_official_gnn_iou}
\end{figure}

The direct single-level DNF reaches \(0.533\) mean mIoU, compared with
\(0.623\) for the complete multilevel pipeline.  This gap remains even though
the single-level baseline receives a larger 10,000-step Adam budget,
with more nominal optimizer updates than the 1200 coarse-grid Adam
updates, eight corrector L-BFGS iterations, and 60 outer mesh-L-BFGS updates in
the complete pipeline.  Direct one-level optimization does not match the multilevel construction under this budget; an Adam and an L-BFGS update differ in cost, so equal-step count is not a wall-clock comparison.  The evaluated PINN
reaches three-class mIoU values of \(0.201\) and \(0.159\) on Data~1 and
Data~2, but does not produce a valid three-class Otsu split on Data~3 or
Data~4; an aggregate PINN mean is therefore not reported.  Relative to the official method, the complete pipeline
improves all four cases, with the largest absolute gain occurring on Data~1.
This comparison evaluates the complete KTC2023 pipeline; the matched LUNDIsim
experiment in Section~\ref{subsec:transfer_ablation} isolates the causal effect
of the ResiPhy-GNN transfer.  The official and DNF results use
the official two-threshold Otsu protocol; as noted above, the PINN fails this
three-class segmentation on Data~3 and Data~4.  Ground-truth labels are used only for evaluation,
not for computing \(\lambda_\sigma\) or fitting the corrector.
\section{Conclusion}
\label{sec:conclusion}

ResiPhy-MDNF combines discrete neural field optimization with a residual-based
GNN transfer operator across nested multilevel spaces.  It uses the prescribed
numerical operator directly, without surrogate models or offline training data.

On Darcy, the two-level path reduces coefficient and state errors by
\(6.76\times\) and \(10.4\times\) over the direct solver at \(25\%\) additional
grid work, and the four-level LUNDIsim pipeline recovers a high-contrast
geological permeability field.  A matched ablation isolates the ResiPhy-GNN
contribution at \(8\)--\(9\%\) beyond interpolation.  On KTC2023 EIT, the
pipeline achieves \(0.623\) mean Otsu mIoU versus \(0.603\) for the official
CEM reconstruction, with gains on all four phantoms.  The same framework
applies to both Darcy and EIT despite differing PDEs, discretizations, and
observation geometries.

The present experiments use regular Cartesian grids.  While the graph-based
transfer does not assume structured connectivity, its effectiveness on
unstructured and adaptively refined meshes remains to be evaluated; extending
the framework to such discretizations is an important direction for future
work. 
\bibliographystyle{cas-model2-names}
\bibliography{references}

\end{document}